\documentclass{article}
\usepackage{amssymb,latexsym,pstricks
,pst-node,pst-coil,pst-plot,
}

\newtheorem{theorem}{Theorem}[section]
\newtheorem{proposition}[theorem]{Proposition}
\newtheorem{lemma}[theorem]{Lemma}

\newtheorem{definition}[theorem]{Definition}

\def\Z{\mathbb{Z}}
\def\R{\mathbb{R}}
\def\C{\mathbb{C}}

\def\CP{\mathbb{CP}}
\def\({\left(}
\def\){\right)}
\def\[{\left[}
\def\]{\right]}

\def\<{\langle}
\def\>{\rangle}
\def\H{{\mathbb H}}
\def\1{{\bf 1}}
\def\i{{\bf i}}
\def\j{{\bf j}}
\def\k{{\bf k}}

\def\proof{\smallskip\noindent {\it Proof --- \ }}

\title{Symmetries of the hypergeometric function $\phantom{}_mF_{m-1}$}

\author{%
Oleg Gleizer\\[.05in]
{\normalsize ogleizer@mac.com}\\[.2in]}

\begin{document}
\maketitle

\begin{abstract}
In this paper, we show that the generalized hypergeometric function $\phantom{}_mF_{m-1}$ 
has a one parameter group of local symmetries, which is a conjugation of a flow of a rational 
Calogero-Mozer system. We use the symmetry to construct fermionic fields on a complex torus, 
which have linear-algebraic properties similar to those of the local solutions of the generalized 
hypergeometric equation. The fields admit a non-trivial action of the quaternions based on 
the above symmetry. We use the similarity between the linear-algebraic structures to introduce 
the quaternionic action on the direct sum of the space of solutions of the generalized 
hypergeometric equation and its dual. As a side product, we construct a ``good'' basis for 
the monodromy operators of the generalized hypergeometric equation inspired by the study 
of multiple flag varieties with finitely many orbits of the diagonal action of the general linear group 
by Magyar, Weyman, and Zelevinsky. As an example of computational effectiveness of the basis, 
we give a proof of the existence of the monodromy invariant hermitian form on the space of 
solutions of the generalized hypergeometric equation (in the case of real local exponents) 
different from the proofs of Beukers and Heckman and of Haraoka. As another side product, 
we prove an elliptic generalization of Cauchy identity.    
\end{abstract} 

\tableofcontents 

%%%%%%%%%%%%%%%%%%%%%%%%%%%%%%%%%%%%%%%%%%%%%%%%%
%%%%%%%%%%%%%%%%%%%%%%%%%%%%%%%%%%%%%%%%%%%%%%%%%
\section{Introduction}
\label{sec:intro}
This section is a short seminar-style exposition of the paper from history and motivations 
to main results to open questions. It should suffice the reader who wants to understand 
the results without going into too many details. \\
  
	One of the ways to define the {\it generalized hypergeometric function} 
$\phantom{}_mF_{m-1}$ is by means of power series:
\begin{equation} \label{equation:HG:F:B-H} 
\phantom{|}_mF_{m-1}\(\left.
\displaystyle{{b_1,\cdots ,b_m} \atop{c_1,\cdots ,c_{m-1}}} \right|\, z\) = 
\sum\limits_{n=0}^{\infty} \displaystyle{\frac{(b_1)_n\cdots (b_m)_n}
{(c_1)_n \cdots (c_{m-1})_n\, n!}}\, z^n.
\end{equation}

\noindent Here $(x)_n$ is the {\it Pochhammer symbol} $(x)_n = \Gamma(x + n)/ \Gamma(x)$. 
The right hand side of (\ref{equation:HG:F:B-H}) converges inside the unit circle of 
the complex plane centered at zero. \\

	The following ODE is known as the {\it generalized hypergeometric equation} 
(GHGE):
\begin{equation} \label{equation:GHGE} 
     z \(z\, \displaystyle{\frac{d}{d\, z}} - c_1\) \cdots 
     \(z\, \displaystyle{\frac{d}{d\, z}} - c_m\)\, f(z) = 
     \(z\, \displaystyle{\frac{d}{d\, z}} - b_1\) \cdots 
     \(z\, \displaystyle{\frac{d}{d\, z}} - b_m\) f(z).
\end{equation}
\medskip

\noindent 	The functions 
\begin{equation} \label{equation:GHGE:sols:0}
  z^{b_j}\! \phantom{|}_mF_{m-1}\(\left. \displaystyle{{b_j-c_1, \cdots , b_j-c_m}
  \atop{b_j-b_1+1,\cdots ,\widehat{b_j-b_j+1},\cdots ,b_j-b_m+1}} \right|\, z\)
\end{equation}
and
\begin{equation} \label{equation:GHGE:sols:infty}
 z^{c_i}\! \phantom{|}_mF_{m-1}\(\left. \displaystyle{{b_1 - c_i, \cdots , b_m - c_i}
 \atop{c_1 - c_i + 1,\cdots ,\widehat{c_i - c_i + 1},\cdots ,c_m - c_i + 1}} \right|\, 
 \displaystyle{\frac1z}\) 
\end{equation}
\medskip

\noindent form bases of the spaces of solutions of the GHGE in a vicinity of zero and 
of infinity respectively, if the {\it local exponents} at zero $b_i$ are distinct $mod~1$ 
as well as the local exponents at infinity $-c_i$. \\

	The {\it monodromy} of the GHGE was found by K.~Okubo in \cite{O} and 
independently by F.~Beukers and G.~Heckman in \cite{BH}. The monodromy 
matrices of the GHGE give probably the most important example of a {\it rigid local system} 
(see page \pageref{def:rigid} for the definition). An algorithm to construct all rigid local systems 
on the Riemann sphere was presented by N.~Katz in \cite{NKatz} and translated into the 
language of linear algebra by M.~Dettweiler and S.~Reiter in \cite{DR}. Y.~Haraoka and T.~Yokoyama in \cite{HY1} give an algorithm to construct all {\it semisimple} rigid local 
systems in the {\it Okubo normal form} and prove that the corresponding Fuchsian systems 
have integral solutions (see also \cite{H1} and \cite{H2}). However, both algorithms are so 
computationally complicated that one should choose a Fuchsian system to apply them or any 
other means of studies to very carefully. Here is one possible criterion to pick Fuchsian systems 
for detailed studies. \\ 

	The local exponents of a Fuchsian system stratify its space of solution into a flag near 
each singularity. There is no basis in the space of solutions which is simultaneously "good" 
for all the flags, so the flags should be considered up to a basis change and thus give rise to 
flag varieties. One way to find  the most important Fuchsian systems is to look for the simplest nontrivial multiple flag varieties. P.~Magyar, J.~Weyman, and A.~Zelevinsky classified in 
\cite{MWZ} all indecomposable multiple flag varieties with finitely many orbits under 
the diagonal action of the general linear group (of simultaneous base changes). 
It turned out that there were three infinite series: the hypergeometric, the odd, and the even 
and two extra cases $E_8$ and ${\hat E}_8$. The Fuchsian systems corresponding to 
all the cases were constructed by the author in \cite{G1}. It is no coincidence that when 
the first major breakthrough in understanding rigid local systems had been made earlier 
by C.~Simpson in \cite{S}, the local systems he had constructed were the hypergeometric, 
the odd, the even, and the extra case ${\hat E}_8$. The same results were obtained 
with different techniques by V.~Kostov in \cite{Ko1}. (See also Kostov's survey of the
{\it Deligne-Simpson problem} in \cite{Ko2}). So in a sense, the Fuchsian systems constructed 
in \cite{G1} are the "more equal animals" from the bestiary of Fuchsian systems on the 
Riemann sphere. Among them, the Fuchsian system corresponding to the hypergeometric 
case is definitely the "most equal animal". On the one hand, it is arguably the simplest 
non-trivial Fuchsian system. On the other hand, it is equivalent to the GHGE as a flat 
connection: it has the same singularities and the same monodromy. We shall call this system 
the $m$-{\it hypergeometric system} (mHGS). It is the main object of study in this paper. 
Any information regarding the mHGS we obtain can (and will) be translated into the 
information regarding the GHGE. 

\begin{definition} \label{def:hgtriple:add}
  Linear operators $A$, $B$, and $C$ acting on a complex linear space
  $\C^m$ are called an \underline{additive hypergeometric triple}, if \\ \\
  $\bullet$~~~$A + B + C = 0$; \\
  $\bullet$~~~$rank\, B = rank\, C = m$; \\
  $\bullet$~~~$A$ is a diagonalizable operator with two different eigenvalues: $a_1$ 
  of multiplicity $1$ and \\
  $\phantom{\bullet~~}$ $a_2$ of multiplicity $m-1$; \\
  $\bullet$~~~the operators are generic within the above restrictions. 
\end{definition} 

	Let $A$, $B$, and $C$ be an additive hypergeometric triple. Let the eigenvalues 
of $B$ and $C$ be $b_1, \dots, b_m$ and $-c_1, \dots ,-c_m$ respectively and let 
$v_1, \dots ,v_m$ and $w_1, \dots ,w_m$ be the corresponding eigenvectors. 
A consequence of the fact that the operators are generic is that  their 
eigenvalues are generic complex numbers in the plane given by the {\it trace condition}:
\begin{equation} \label{equation:trace}
  a_1 + (m-1) a_2 + \sum\limits_{i=1}^m b_i - c_ i =0.
\end{equation}

\noindent In particular, none of the differences $b_i - b_j$ and $c_i - c_j$ is an integer 
for $i \ne j$, and neither is $b_i - c_j$ for all $i$ and $j$. \\

	Let $u$ be the eigenvector of $A$ corresponding to the eigenvalue $a_1$. 
Another consequence of the fact that that the operators are generic is that there exists 
a unique way to choose $v_i$ and $w_i$ so that
$$ \sum\limits_{i = 1}^m v_i = \sum\limits_{i = 1}^m w_i = u. $$ \medskip
 
	Let $f: \CP^1 \setminus\{0,1,\infty\} \to \C^m$. The mHGS is the following 
Fuchsian system:
\begin{equation} \label{equation:m-hgs}
  \frac{d\,f}{d\,z}=\[\frac{A}{z-1} + \frac{B}{z}\]f(z).
\end{equation}

\noindent The additive hypergeometric triple is the triple of the {\it residue operators} 
of the mHGS. $C$ is the residue at infinity, so it is not explicitly visible in 
(\ref{equation:m-hgs}). We shall call the space where the residues act the 
{\it residue space}. The following feature of the residue space was proven in \cite{G1}:
\begin{theorem} \label{thm:(v,v),(w,w)}
  There exists a unique up to a constant multiple complex symmetric scalar 
  product $(*,*)_r$ on the residue space of the mHGS such that the bases 
  $v_i$ and $w_i$ are simultaneously orthogonal with respect to it, given by
  \begin{equation} \label{equation:(v,v),(w,w)}
    \(v_i, v_j\)_r = \delta_{ij}\, \nu_{i\, r}^{2}~~~ \mbox{and}~~~ 
    \(w_i, w_j\)_r = \delta_{ij}\, \mu_{i\, r}^{2},
  \end{equation}

  \noindent where  
  \begin{equation} \label{equation:mu_0,nu_0}
      \mu_{i\, r}^2 = \displaystyle{\frac{\prod\limits_{k=1}^m 
      (a_2 + b_k - c_i)}{\prod\limits_{k=1\atop{k \ne i}}^m (c_k-c_i)}}~~~ \mbox{and}~~~  
      \nu_{i\, r}^2 = \displaystyle{\frac{\prod\limits_{k=1}^m (a_2 + b_i - c_k)}
      {\prod\limits_{k=1\atop{k \ne i}}^m (b_i-b_k)}}.
  \end{equation}
\end{theorem}

	When all the local exponents are real, the form $(*,*)_r$ is 
real symmetric. It is important to know when the real form is sign-definite. 
The following lemma is proven in \cite{G1}:

\begin{lemma} \label{lemma:hg:form:pos}
  Renumbering if necessary, we can think that $b_1 < \dots < b_m$ and 
  $c_1 < \dots < c_m$. If the form $\epsilon\, (*,*)_r$ is positive-definite 
  for $\epsilon=\pm 1$, then $\epsilon=sign(a_2 - a_1)$. If the inequalities of the first column 
  hold, then $a_2 > a_1$. If the inequalities of the second column hold, then $a_2 < a_1$.
  \begin{equation} \label{equation:hg:form:pos}
    \begin{array}{|c|c|}
      \hline
       & \\
      \begin{array}{lcccl}
        c_1 - b_1              & < & a_2    & < & c_2 - b_1        \\
        c_2 - b_2              & < & a_2    & < & c_3 - b_2        \\
        \vdots                    &   & \vdots &   & \vdots               \\
        c_{m-1} - b_{m-1} & < & a_2    & < & c_m - b_{m-1} \\
        c_m - b_m            & < & a_2
      \end{array} &

      \begin{array}{lcccl}
        c_{m-1} - b_m       & < & a_2    & < & c_m - b_m     \\
        c_{m-2} - b_{m-1} & < & a_2    & < & c_{m-1} - b_{m-1} \\
        \vdots                   &     & \vdots &   & \vdots      \\
        c_1 - b_2              & <  & a_2    & < & c_2  - b_2 \\
        \phantom{c_{m-2} -  b_{m-1}} &   & a_2 & < & c_1 - b_1
      \end{array} \\
      & \\
      \hline
    \end{array}.
  \end{equation}
\end{lemma}
\medskip

\noindent Note that we can extend the sign-definite real symmetric form $(*,*)_r$
to complex numbers either in the complex symmetric or in the hermitian way. \\

	Let $V$ be the operator of the basis switch from the standard basis to the basis 
$v_i$. Let
\begin{equation} \label{equation:S} 
  S = V X V^{-1}, 
\end{equation}

\noindent where $X$ is the following $m \times m$ matrix:
\begin{equation} \label{equation:X}
  X_{ij} = \left\{
 \begin{array}{ccl}
   \displaystyle{\frac{1}{b_j - b_i}}, & if & i \ne j, \\
    & & \\
   -\sum\limits_{j=1 \atop{j \ne i}}^m \displaystyle{\frac{1}{b_j - b_i}}, & if & i=j.
 \end{array} \right.
\end{equation}

\noindent The following two systems with complex times play a major role in our investigation:
\begin{equation} \label{equation:main}
  \begin{array}{lll}
      {\dot B}(\tau_1) + \[B(\tau_1),S\] = k_1 Id, & \phantom{hmm} & 
      {\dot C}(\tau_2) - \[C(\tau_2),S\] = k_2 Id, \\
       & & \\
      B(0) = B, & & C(0) = C.
   \end{array}
\end{equation}

\noindent Here $[*,*]$ is the usual commutator, $k_1$ and $k_2$ are arbitrary complex 
constants, and $Id$ is the $m \times m$ identity matrix. The following theorem shows the 
importance of the systems (\ref{equation:main}) for studying the mHGS.

\begin{theorem} \label{thm:Frob_series}
If $a_2 = 0$, then mHGS has ``nice'' Fr\"obenius series solutions 
in a local parameter $z$ near zero and in a local parameter $z^{-1}$ near infinity 
respectively:
\begin{equation} \label{equation:T0}
  (T_0)_i = z^{b_i} \(\sum_{n=0}^{\infty} \alpha_{in}\, e^{S n}\, z^n\) v_i
  ~~~ \mbox{and}
\end{equation}

\begin{equation} \label{equation:Tinfty}
  (T_{\infty})_i = z^{c_i} \(\sum_{n=0}^{\infty} \beta_{in}\, e^{-S n}\, z^{-n}\) w_i,
  ~~~ \mbox{where}
\end{equation}

\begin{equation}
  \label{equation:alpha,beta}
  \begin{array}{ll}
     \alpha_{in} = \prod\limits_{k=1}^m 
     \displaystyle{\frac{\Gamma(b_i - c_k + n)\, \Gamma(b_i - b_k + 1)}
     {\Gamma(b_i - c_k)\, \Gamma(b_i - b_k + n+1)}}, &
     \beta_{in} = \prod\limits_{k=1}^m 
     \displaystyle{\frac{\Gamma(b_k - c_i + n)\, \Gamma(c_k - c_i + 1)}
     {\Gamma(b_k - c_i)\, \Gamma(c_k - c_i + n+1)}}. 
  \end{array}
\end{equation}

\noindent The series (\ref{equation:T0}) and (\ref{equation:Tinfty}) converge 
inside the unit circles centered at zero and at infinity. 
\end{theorem}

\begin{lemma} \label{lemma:ev} 
The eigenvalues of $B(\tau_1)$ and $C(\tau_2)$ evolve linearly: 
$b_i(\tau_1) = b_i + k_1 \tau_1$ and \\ $-c_i(\tau_2) = -c_i + k_2 \tau_2$. 
The corresponding eigenvectors are
  \begin{equation} \label{equation:v(t),w(t)}
    v_i(\tau_1) = e^{S \tau_1} v_i,~~~ w_i(\tau_2) = e^{-S \tau_2} w_i,
  \end{equation}

\noindent where $v_i$ and $w_i$ are the eigenvectors of the original operators $B$ and $C$. 
\end{lemma}

\noindent Thus, up to normalizing factors, the coefficients of the local solutions 
(\ref{equation:T0}) and (\ref{equation:Tinfty}) are the eigenvectors $v_i(\tau_1)$ and 
$w_i(\tau_2)$ of the operators $B(\tau_1)$ and $C(\tau_2)$ for the integral times 
$\tau_1 = n$ and $\tau_2 = -n$. \\

	The GHGE has the subspace of solutions holomorphic at one of dimension $m-1$, 
which corresponds to the case $a_2=0$ for the mHGS. However, from the point of view 
of representation theory, the case of the traceless residue operators seems to be more 
important. Unfortunately, when $a_2 \ne 0$, the nice Fr\"obenius series for the local 
solutions similar to the above exist no more. To study the mHGS in this case, one needs 
an apparatus different from the residue calculus. \\ 

	In the case $a_2 =0$, the study of the series (\ref{equation:T0}) and (\ref{equation:Tinfty}) 
is completely parallel to the study of the local solutions (\ref{equation:GHGE:sols:0}) and 
(\ref{equation:GHGE:sols:infty}) of the GHGE. However, our approach through the dynamical 
systems (\ref{equation:main}) suggests that taking the times $\tau_1$ and $\tau_2$ along two 
linearly dependent vectors in the complex plane is a trigonometric degeneration of a more 
general elliptic situation. So on the one hand, it seems desirable to construct an elliptic analogue 
of (\ref{equation:m-hgs}). On the other hand, we know from \cite{SJM} that a solution to 
(\ref{equation:m-hgs}) can be realized as a fermionic field. We shall use this idea to treat
the case $a_2 \ne 0$ by means of linear algebra rather than the residue calculus. 

\begin{definition} \label{def:hgtriple:mult}
  Linear operators $M_0$, $M_1$, and $M_{\infty}$ acting on a complex linear space
  $\C^m$ are called a \underline{multiplicative hypergeometric triple}, if \\ \\
  $\bullet$~~~$M_{\infty} M_1 M_0 = Id$; \\
  $\bullet$~~~$rank\, M_0 = rank\, M_1 = m$; \\
  $\bullet$~~~$M_1$ is a diagonalizable operator with two different eigenvalues 
  of multiplicities $1$ and $m-1$; \\
  $\bullet$~~~the operators are generic within the above restrictions. 
\end{definition}

	Quite obviously, the monodromy operators of the mHGS form a multiplicative 
hypergeometric triple. The eigenvalues of $M_0$ and $M_{\infty}$ are 
$e^{2 \pi \sqrt{-1}\, b_i}$ and $e^{-2 \pi \sqrt{-1}\, c_i}$ respectively. $M_1$ is 
diagonalizable and has two eigenvalues: $e^{2 \pi \sqrt{-1}\, a_1}$ of multiplicity one 
and $e^{2 \pi \sqrt{-1}\, a_2}$ of multiplicity $m-1$. The triple loop on 
$\CP^1 \setminus\{0, 1, \infty\}$ passing around zero, one, and infinity once in the positive
(counter clockwise) direction is contractible, so $M_{\infty} M_1 M_0 = Id$. Note that the 
product of the operators is taken in the opposite order: that is because they act on the right
taking linear combinations of columns of the fundamental matrix. \\

	Let $p_i$ and $q_i$ be the eigenvectors of $M_0$ and $M_{\infty}$ corresponding 
to the eigenvalues $e^{2 \pi \sqrt{-1}\, b_i}$ and $e^{-2 \pi \sqrt{-1}\, c_i}$ respectively. 
Let $r$ be the eigenvector of $M_1$ corresponding to the eigenvalue $a_1$. The operators 
are generic, so it is uniquely possible to choose the $p_i$ and $q_i$ so that 
\begin{equation} \label{equation:normal:trig}
  \sum\limits_{i = 1}^m p_i = \sum\limits_{i = 1}^m q_i = r.
\end{equation}

	The following theorem is proven in \cite{BH} for the GHGE. We reprove this theorem 
for the mHGS with any $a_2$ in this paper. The basis we work in is different from the 
basis used by Beukers and Heckman and is more convenient in some situations. 

\begin{theorem} \label{thm:(p,p),(q,q)}
  If all the local exponents are real, then there exists a unique up to a constant multiple 
  monodromy invariant hermitian product $(*,*)_{trig}$ on the space of solutions such that 
  the bases $p_i$ and $q_i$ are simultaneously orthogonal with respect to it, given by
  \begin{equation} \label{equation:(p,p),(q,q)}
    (p_i, p_j)_{trig} = \delta_{ij}\, \nu_{i\, trig}^2~~~ \mbox{and}~~~ 
    (q_i, q_j)_{trig} = \delta_{ij}\, \mu_{i\, trig}^2,
  \end{equation}

  \noindent where  
  \begin{equation} \label{equation:mu_ell,nu_e}
    \mu_{i\, trig}^2 = \displaystyle{\frac{\prod\limits_{k=1}^m 
    \sin \pi (a_2 + b_k - c_i)}{\prod\limits_{k=1\atop{k \ne i}}^m \sin \pi (c_k-c_i)}}
    ~~~ \mbox{and}~~~  
    \nu_{i\, trig}^2 = \displaystyle{\frac{\prod\limits_{k=1}^m \sin \pi (a_2 + b_i - c_k)}
    {\prod\limits_{k=1\atop{k \ne i}}^m \sin \pi (b_i-b_k)}}.
  \end{equation}
\end{theorem}

\noindent The initial motivation for the author to study the subject was to understand 
why the formulae (\ref{equation:(v,v),(w,w)}) and (\ref{equation:(p,p),(q,q)}) look so similar 
despite the different nature of the products: the first is complex
symmetric whereas the second is hermitian. \\ 

	The object which explains this similarity is the following: further in the paper we construct 
four $m$-tuples of fermionic fields $\(F_0\)_i$, $\(F_0\)_i^{\dagger}$, $\(F_{\infty}\)_j$, 
and $\(F_{\infty}\)_j^{\dagger}$ such that their {\it vacuum expectation values} 
(see page \pageref{equation:Wick's rule} for the definition) are 
\begin{equation} \label{equation:1/sn}
    \left\< \(F_{\infty}\)_j^{\dagger}\, \(F_0\)_i \right\> = 
    \displaystyle{\frac{1}{sn(a_2 + b_i - c_j)}} =
    \left\< \(F_{\infty}\)_j\, \(F_0\)_i^{\dagger} \right\>.
\end{equation}

\noindent Here $sn$ is the elliptic sine of Jacobi and $a_2$, $b_i$ and $-c_i$ are 
the local exponents of the mHGS. We use the pairing (\ref{equation:1/sn}) to identify 
the linear spaces spanned by $\(F_0\)_i$ and $\(F_{\infty}\)_i^{\dagger}$ as well as by $\(F_0\)_i^{\dagger}$ and $\(F_{\infty}\)_i$. Let us call the resulting spaces $\H$ and $\H'$
respectively. These spaces resemble the residue space and the space of solutions of 
the mHGS in many ways. In particular, if we rewrite the formulae (\ref{equation:(v,v),(w,w)}) 
and (\ref{equation:(p,p),(q,q)}) as
$$
  \begin{array}{lll}
    \(\displaystyle{\frac{v_i}{\nu_{i\, r}^{2}}}, \displaystyle{\frac{v_j}{\nu_{j\, r}^{2}}}\)_r = 
    \displaystyle{\frac{\delta_{ij}}{\nu_{i\, r}^2}} & \phantom{trr-trr} &
    \(\displaystyle{\frac{w_i}{\mu_{i\, r}^{2}}}, \displaystyle{\frac{w_j}{\mu_{j\, r}^{2}}}\)_r = 
    \displaystyle{\frac{\delta_{ij}}{\mu_{i\, r}^2}} \\ & & \\
    \(\displaystyle{\frac{p_i}{\nu_{i\, trig}^{2}}}, 
    \displaystyle{\frac{p_j}{\nu_{j\, trig}^{2}}}\)_{trig} = 
    \displaystyle{\frac{\delta_{ij}}{\nu_{i\, trig}^2}} & &
    \(\displaystyle{\frac{w_i}{\mu_{i\, trig}^{2}}}, 
    \displaystyle{\frac{w_j}{\mu_{j\, trig}^{2}}}\)_{trig} = 
    \displaystyle{\frac{\delta_{ij}}{\mu_{i\, trig}^2}}, 
  \end{array}
$$ 

\noindent then Theorems \ref{thm:(v,v),(w,w)} and \ref{thm:(p,p),(q,q)} start looking
very similar to the following:
\begin{theorem} \label{thm:(F,F):ell}
There exists a unique up to a constant multiple complex symmetric scalar 
product $(*,*)_{ell}$ on $\H$ such that the bases $\(F_0\)_i$ and 
$\(F_{\infty}\)^{\dagger}_i$ are simultaneously orthogonal with respect to it, given by
\begin{equation} \label{equation:(F,F):ell}
  \(\(F_0\)_i, \(F_0\)_j\)_{ell} = \displaystyle{\frac{\delta_{ij}}{\nu_{i\, ell}^2}}~~~ 
  \mbox{and}~~~ 
  \(\(F_{\infty}\)_i^{\dagger}, \(F_{\infty}\)_j^{\dagger}\)_{ell} = 
  \displaystyle{\frac{\delta_{ij}}{\mu_{i\, ell}^2}},
\end{equation}

\noindent where  
\begin{equation} \label{equation:mu_ell,nu_ell}
  \mu_{i\, ell}^2 = \displaystyle{\frac{\prod\limits_{k=1}^m 
  sn(a_2 + b_k - c_i)}{\prod\limits_{k=1\atop{k \ne i}}^m sn(c_k-c_i)}}~~~ \mbox{and}~~~  
  \nu_{i\, ell}^2 = \displaystyle{\frac{\prod\limits_{k=1}^m sn (a_2 + b_i - c_k)}
  {\prod\limits_{k=1\atop{k \ne i}}^m sn (b_i-b_k)}},
\end{equation}

\noindent Similar formulae hold for $\H'$.
\end{theorem}

	The eigenvalues of the monodromy operators of a linear regular system are more 
important than the local exponents due to the following reason. Consider the linear regular 
matrix system 
\begin{equation} \label{equation:gauge}
  \displaystyle{\frac{dT(z)}{dz}} = R(z) T(z).
\end{equation} 
The {\it gauge transformation} $T \mapsto g(z) T$ replaces $R$ by
$$
  \displaystyle{\frac{dg}{dz}} g^{-1} + gRg^{-1}.
$$
Most often $g$ are taken holomorphic and holomorphically invertible away from the poles
of the original system, so that the new system has the same singularities as the old one. 
The only invariant under such a transformation is the monodromy group of the system; 
the residue matrices are not preserved. Combining this perspective with the real local exponents, 
we can think that the local exponents belong to the semi-interval $[0,1)$ right away. 
Renumbering if necessary, we can, similarly to Lemma \ref{lemma:hg:form:pos}, think that 
$0 \le b_1 < \dots < b_m < 1$ and $0 \le c_1 < \dots < c_m < 1$.\\

	Let the period lattice of the elliptic sine in (\ref{equation:1/sn}) be generated by
$2 \omega_1$ and $\omega_2$, see page \pageref{equation:normSn}. The fact that 
unlike the residue space and the space of solutions, the spaces $\H$ and $\H'$ 
naturally come in a pair, gives rise to the following two parts theorem:
\begin{theorem} \label{thm:hk} $\phantom{trr}$ \\ 
 1.~~ Let us restrict ourselves to the case $\omega_1 = 1$ and $\omega_2 = \sqrt{-1} s$, 
  where $s \in \R$. Let $a_2$, $b_i$ and $c_i$ be generic real numbers from $[0,1)$ 
  satisfying the positivity conditions (\ref{equation:hg:form:pos}). Then the formulae 
  (\ref{equation:(F,F):ell}) give hermitian forms on $\H$ and $\H'$. If the inequalities 
  (\ref{equation:hg:form:pos}) are not satisfied, then the forms have nontrivial signatures. \\ 
  
  \noindent 2.~~ The quaternions act on $\H \oplus \H'$ by means of the the following formulae:
  \begin{equation}\label{equation:quaternions:F}
    \begin{array}{lll}
      \i \(F_0\)_i =  -\(F_0\)_i^{\dagger} &
      \j \(F_0\)_i = -\sqrt{-1}\, \(F_0\)_i & 
      \k \(F_0\)_i = \sqrt{-1}\, \(F_0\)_i^{\dagger} \\
       & & \\
      \i \(F_0\)_i^{\dagger} = \(F_0\)_i & 
      \j \(F_0\)_i^{\dagger} = \sqrt{-1}\, \(F_0\)_i^{\dagger} & 
      \k \(F_0\)_i^{\dagger} = \sqrt{-1}\, \(F_0\)_i \\
       & & \\ 
      \i \(F_{\infty}\)_i = -\(F_{\infty}\)_i^{\dagger} & 
      \j \(F_{\infty}\)_i = -\sqrt{-1}\, \(F_{\infty}\)_i & 
      \k \(F_{\infty}\)_i = \sqrt{-1}\, \(F_{\infty}\)_i^{\dagger} \\
       & & \\
      \i \(F_{\infty}\)_i^{\dagger} = \(F_{\infty}\)_i & 
      \j \(F_{\infty}\)_i^{\dagger} = \sqrt{-1}\, \(F_{\infty}\)_i^{\dagger} &  
      \k \(F_{\infty}\)_i^{\dagger} = \sqrt{-1}\, \(F_{\infty}\)_i   
    \end{array}
  \end{equation}
  \smallskip
  
   \noindent If the hermitian forms are sign-definite, then the quaternionic action is 
   hyperk\"ahler, which explains the simultaneous presence of the complex
   symmetric and hermitian forms, see \cite{Hit1}. 
  \end{theorem}
\smallskip

	A construction of the quaternionic action (\ref{equation:quaternions:F}) further 
in the paper shows that the spaces $\H$ and $\H'$ are naturally dual to each other. 
Moreover, the sum of the solution space and its dual turns out to be equivalent to the 
trigonometric limit of $\H \oplus \H'$, whereas the sum of the residue space and its dual 
is equivalent to the rational limit. The following commuting diagram illustrates how 
it works: 
$$ \label{equation:lim_diag}
  \begin{array}{ccccc}
    \mbox{?} & \simeq & \left\< \(F_{\infty}\)_j^{\dagger}\, \(F_0\)_i \right\> & = &
    \displaystyle{\frac{1}{sn(a_2 + b_i - c_j)}} \\ & & & & \\ 
     & & trig \downarrow & \omega_2 \to \infty & \downarrow trig \\ & & & & \\
    \(\displaystyle{\frac{q_j}{\mu^2_{j\, trig}}},\displaystyle{\frac{p_i}{\nu^2_{i\, trig}}}\)_{trig} & 
    \simeq & \left\< \(F_{\infty}\)_{j\, trig}^{\dagger}\, \(F_0\)_{i\, trig} \right\> & = & 
    \displaystyle{\frac{\pi}{\sin \pi (a_2 + b_i - c_j)}} \\ & & & & \\ 
     & & r \downarrow & \omega_1 \to \infty & \downarrow r \\ & & & & \\ 
    \(\displaystyle{\frac{w_j}{\mu^2_{j\, r}}},\displaystyle{\frac{v_i}{\nu^2_{i\, r}}}\)_r &  
    \simeq & \(\displaystyle{\frac{1}{\mu_{j\, r}^2}}\, 
    \displaystyle{\frac{w_j \oplus -\sqrt{-1}\, w_j}{\sqrt{2}}}, 
    \displaystyle{\frac{1}{\nu_{i\, r}^2}}\, 
    \displaystyle{\frac{v_i \oplus \sqrt{-1}\, v_i}{\sqrt{2}}}\) & = & 
    \displaystyle{\frac{1}{a_2 + b_i - c_j}} 
  \end{array}
$$
\smallskip

\noindent The equivalences in the diagram occur at the level of formal linear algebra. 
Geometrically, $p_i$ and $q_j$ are multivalued functions on $\CP^1 \setminus \{0, 1, \infty\}$, 
which can be expressed by convergent Fr\"obenius series or by Euler integrals, undergo 
monodromy transformations, etc. Their counterparts $\(F_0\)_{i\, trig}$ and 
$\(F_{\infty}\)_{j\, trig}^{\dagger}$ are infinite sums of annihilation/creation operators, which 
obtain geometric meaning when applied to vectors of the corresponding Fock space. 
The vectors $v_i$ and $w_j$ are the eigenvectors of the residue operators $B$ and $C$ of 
the mHGS, whereas the $v_i \oplus \sqrt{-1}\, v_i$ and $w_j \oplus -\sqrt{-1}\, w_j$ should be 
thought of as isotropic generators of some Clifford algebra. Still, some linear-algebraic 
information can be dragged from one side to the other in a meaningful fashion. In particular, 
the quaternionic action (\ref{equation:quaternions:F}) can be introduced on the sums of 
the solution and the residue spaces of the mHGS with their respective duals, explaining that 
in fact both the complex symmetric and the hermitian structures are present in both cases. 
The limiting procedure explains the similarity between the formulae (\ref{equation:(v,v),(w,w)}) 
and (\ref{equation:(p,p),(q,q)}) for the complex symmetric product $(*,*)_r$ on the residue space 
and the hermitian product $(*,*)_{trig}$ on the space of solutions of the mHGS. \\

	The question mark in the diagram stands for the unknown Fuchsian system on the torus, 
which becomes the mHGS as the period $\omega_2 \to \infty$. A construction 
of such a system would deepen the above equivalences from the linear-algebraic to the 
geometric level. In particular, it should allow to drag the action of the monodromy group 
of the mHGS to the field side of the diagram in a meaningful geometric fashion. The author plans 
to construct this system in a subsequent publication.

\section{More results} %%%%%%%%%%%%%%%%%%%%%%%%%%%%%%%%%%%%%
%%%%%%%%%%%%%%%%%%%%%%%%%%%%%%%%%%%%%%%%%%%%%%%%
In this section, some further results are presented. We also give a couple of proofs, 
which are essential for understanding the introduced objects and some of the results. 

\subsection{Symmetries of the mHGS and the Calogero - Moser flow} %%%%%%%%%%%%%
For reader's convenience, some basic information about the rational Calogero-Mozer system 
is provided in the Appendix on page \pageref{sec:CM}. \\ 

	The local solutions (\ref{equation:T0}) and (\ref{equation:Tinfty}) of the mHGS at zero 
and at infinity are constructed in terms of the operator $S = V X V^{-1}$. The matrix $X$ 
(\ref{equation:X}) depends on the eigenvalues $b_i$ of the residue operator $B$ only. $V$ is 
the operator of the basis switch from the standard basis to the basis $v_i$ of the eigenvectors 
of $B$. To restore the symmetry of the construction, let us consider the following $m \times m$ 
matrix $Y$: 
\begin{equation} \label{equation:Y}
Y_{ij} = \left\{
 \begin{array}{ccl}
   \displaystyle{\frac{1}{c_i - c_j}}, & if & i \ne j, \\
    & & \\
   -\sum\limits_{j=1 \atop{j \ne i}}^m \displaystyle{\frac{1}{c_i - c_j}}, & if & i=j.
 \end{array} \right.
\end{equation}

\begin{lemma} \label{lemma:propXY} $\phantom{}$ \\ \\ 
  $\bullet$~~ The Jordan normal form of both $X$ and $Y$ is a single block of full size 
  $m$ with the eigenvalue zero. The corresponding eigenvector of both $X$ and $Y$ is 
  $e = (1, \dots, 1)$. \\
  
  \noindent $\bullet$~~ Let $B_d = diag(b_1, \dots ,b_m)$ and $C_d = diag(-c_1, \dots ,-c_m)$ 
  be the Jordan normal forms of the operators $B$ and $C$ respectively. Then the matrices 
  $X$ and $Y$ satisfy a special type of the rational Calogero-Moser equation: 
   \begin{equation} \label{equation:C-M}
          \[X, B_d\] = \[Y, C_d\] = e \otimes e - Id.
   \end{equation}

  \noindent $\bullet$~~ Let $W$ be the operator of the basis switch from the standard basis to the 
  basis $w_i$ of the eigenvectors of the residue operator $C$. Then  
  \begin{equation}
    \label{equation:VXYW}
    V X V^{-1} = -W Y W^{-1}.
  \end{equation} 
\end{lemma}
\medskip

\noindent The last formula of the lemma shows that $X$ and $-Y$ are matrices of the operator 
$S$ from the dynamical equations (\ref{equation:main}) in the bases $v_i$ and $w_i$ respectively. 
 \label{S} The formula $X e = 0$ is equivalent to $S u = 0$, which implies that 
\begin{equation} \label{equation:normal:r}
  \sum\limits_{i = 1}^m v_i(\tau_1) = \sum\limits_{i = 1}^m w_i(\tau_2) = u.
\end{equation}

\noindent for any complex times $\tau_1$ and $\tau_2$. \\

	Recall that the vectors $v_i(\tau_1)$ and $w_i(\tau_2)$ are the eigenvectors of the 
operators $B(\tau_1)$ and $C(\tau_2)$ corresponding to the eigenvalues $b_i(\tau_1)$ 
and $-c_i(\tau_2)$ respectively, see Lemma \ref{lemma:ev}. Let us introduce the following 
notations:
\begin{equation}
\label{equation:mu,nu}
\mu_i^2(\tau) = \displaystyle{\frac{\prod\limits_{k=1}^m 
(a_2 + b_k - c_i + \tau)}{\prod\limits_{k=1\atop{k \ne i}}^m (c_k-c_i)}},~~~~~ 
\nu_i^2(\tau) = \displaystyle{\frac{\prod\limits_{k=1}^m (a_2 + b_i - c_k + \tau)}
{\prod\limits_{k=1\atop{k \ne i}}^m (b_i-b_k)}}.
\end{equation}

\noindent The following theorem looks very similar to Theorems \ref{thm:(v,v),(w,w)}, 
\ref{thm:(p,p),(q,q)}, and \ref{thm:(F,F):ell}. In fact, it underlies all of them:
\begin{theorem} \label{thm:product}
  For any  complex times $\tau_1$ and $\tau_2$, there exists a unique up to 
  a constant multiple complex symmetric scalar product $(*,*)_{\tau_1, \tau_2}$ 
  on $\C^m$ such that the bases $v_i(\tau_1)$ and $w_i(\tau_2)$ 
  are simultaneously orthogonal with respect to it, given by
  \begin{equation} \label{equation:(*,*)_tau12} 
    \begin{array}{l}
       \(v_i(\tau_1), v_j(\tau_1)\)_{\tau_1, \tau_2} = 
       \delta_{ij}\, \nu_i^2\(\tau_1 + \tau_2\), \\ \\
       \(w_i(\tau_2), w_j(\tau_2)\)_{\tau_1, \tau_2} = 
       \delta_{ij}\, \mu_i^2\(\tau_1 + \tau_2\).  
    \end{array} 
  \end{equation}
  
\noindent Also,
\begin{equation}
  \label{equation:(v,w)}
  \(v_i(\tau_1), w_j(\tau_2)\)_{\tau_1, \tau_2} = 
  \displaystyle{\frac{\nu_i^2\(\tau_1 + \tau_2\) \mu_j^2\(\tau_1 + \tau_2\)}
  {a_2 + b_i - c_j + \tau_1 + \tau_2}}.
\end{equation}
\end{theorem}
\bigskip

	Let us define the operator $A(\tau_1, \tau_2)$ by the formula 
$A(\tau_1, \tau_2) + B(\tau_1) + C(\tau_2) = 0$.

\begin{lemma} \label{lemma:A} $\phantom{trr}$ \\ \\ 
1.~ $A(\tau_1, \tau_2)$ is a diagonalizable operator with the eigenvalue 
\begin{equation}
  \label{equation:a_1}
  a_1(\tau_1, \tau_2) = a_1 + (1 - k_1 - m) \tau_1 + (1 - k_2 - m) \tau_2  
\end{equation} 
of multiplicity $1$ and the eigenvalue
\begin{equation}
  \label{equation:a_2}
  a_2(\tau_1, \tau_2) = a_2 + (1 - k_1) \tau_1 + (1 - k_2) \tau_2
\end{equation} 
of multiplicity $m-1$. \\

\noindent 2.~ The vector $u$ is the eigenvector of $A(\tau_1, \tau_2)$ corresponding to 
the eigenvalue $a_1(\tau_1, \tau_2)$. \\ 

\noindent 3.~ For any vector $x \in \C^m$, 
\begin{equation} \label{equation:Ax}
  A(\tau_1, \tau_2)\, x = a_2(\tau_1, \tau_2)\, x - (x,u)_{\tau_1, \tau_2} u
\end{equation}

\noindent 4.~ $(u, u)_{\tau_1, \tau_2} = a_2(\tau_1, \tau_2) - a_1(\tau_1, \tau_2)$. 
\end{lemma}
\bigskip

	Let $H_{\tau_1, \tau_2}$ be a copy of $\C^m$ with two distinguished bases:
$v_i(\tau_1)$ and $w_i(\tau_2)$, equipped with the complex symmetric scalar 
product (\ref{equation:(*,*)_tau12}). The space $H_{\tau_1, \tau_2} \oplus H_{\tau_2, \tau_1}$ 
admits the natural action of the quaternions: 
\begin{equation} \label{equation:quaternions}
  \begin{array}{l}
    \i = \[\begin{array}{cc} 0 & \sqrt{-1}\, e^{S(\tau_1 - \tau_2)} \\
    \sqrt{-1}\, e^{S(\tau_2 - \tau_1)} & 0 \end{array}\], \\ \\
    \j = \[\begin{array}{cc} 0 & -e^{S(\tau_1 - \tau_2)} \\
    e^{S(\tau_2 - \tau_1)} & 0 \end{array}\], \\ \\
    \k = \[\begin{array}{cc} \sqrt{-1} & 0 \\ 0 & -\sqrt{-1} \end{array}\].
  \end{array}
\end{equation}

\noindent The quaternions act on $H_{-\tau_2, -\tau_1} \oplus H_{-\tau_1, -\tau_2}$ 
by means of (\ref{equation:quaternions}) as well. The action on 
$H_{\tau_1, \tau_2} \oplus H_{\tau_2, \tau_1}$ preserves the scalar 
product $(*,*)_{\tau_1,\tau_2} \oplus (*,*)_{\tau_2,\tau_1}$ up to the sign. 
Over the reals, the space $H_{\tau_1, \tau_2} \oplus H_{\tau_2, \tau_1}$ 
is spanned by $v_i(\tau_1)$, $v_i(\tau_2)$, $\sqrt{-1}\, v_i(\tau_1)$, and 
$\sqrt{-1}\, v_i(\tau_2)$. So in the case when the form $(*,*)_{\tau_1,\tau_2}$ 
(and thus $(*,*)_{\tau_2,\tau_1}$) is sign-definite, the quaternionic action 
is hyperk\"ahler. The formulae (\ref{equation:quaternions}) also show that the spaces 
$H_{\tau_1, \tau_2}$ and $H_{\tau_2, \tau_1}$ are dual to each other: the exponents 
acting on them have opposite signs. 

\subsection{The GHGE, the mHGS and the elliptic Cauchy identity} %%%%%%%%%%%%%%

	Since we are going to take a closer look at the GHGE, let us temporarily restrict ourselves 
to the case $a_2 = 0$ when working with the mHGS. \\

	It is not clear from \cite{BH} how the vectors $p_i$ and $q_i$ of 
Theorem \ref{thm:(p,p),(q,q)} correspond to the solutions (\ref{equation:GHGE:sols:0}) 
and (\ref{equation:GHGE:sols:infty}) of the GHGE at zero and at infinity. The following 
lemma relates the classical analytic approach to the linear algebra of \cite{BH}: 
\begin{lemma} \label{lemma:p<->q}
  $$
    p_j = \sum\limits_{i=1}^m \displaystyle{\frac{e^{\pi \sqrt{-1}\, (a_1 + b_j - c_i)}\, 
    \nu_{j\, trig}^2}{\sin \pi (a_2 + b_j - c_i)}}\, q_i
  $$
\end{lemma}
\medskip

	The formula:
\begin{equation} \label{equation:analytic_cont_of_mFm-1_to_infty}
  \begin{array}{l}
    \phantom{}_mF_{m-1}\(\left. \displaystyle{{b_1,\dots ,b_m} \atop{c_1,\dots ,c_{m-1}}} 
    \right| z\) = \sum\limits_{i=1}^m \prod\limits_{k=1\atop{k \ne i}}^m 
    \displaystyle{\frac{\Gamma(b_k - b_i)}{\Gamma(b_k)}}\,
    \prod\limits_{k=1}^{m-1} \displaystyle{\frac{\Gamma(c_k)}{\Gamma(c_ k - b_i)}}\, \times \\ \\
    (-z)^{-b_i}\,\phantom{}_mF_{m-1}\(\left.
    \displaystyle{{b_i - c_1 + 1,\dots ,b_i - c_{m-1} + 1,\, b_i}
    \atop{b_i - b_1 + 1,\dots,\widehat{b_i - b_i + 1},
    \dots ,b_i - b_m + 1}} \right| \displaystyle{\frac{1}{z}} \) 
  \end{array}
\end{equation}

\noindent is proven in \cite{WW} for $m = 2$. Their argument works for any $m$ without 
major change. Using the famous reflection formula
\begin{equation} \label{equation:reflection:gamma}
  \Gamma(z) \Gamma(1-z) = \displaystyle{\frac{\pi}{\sin \pi z}},
\end{equation}
we rewrite (\ref{equation:analytic_cont_of_mFm-1_to_infty}) as 
\begin{equation} \label{equation:g.h.g.eq.:an.cont.}
  \begin{array}{l}
    z^{b_j}\! \phantom{|}_mF_{m-1}\(\left. \displaystyle{{b_j-c_1, \cdots , b_j-c_m}
    \atop{b_j-b_1+1,\cdots ,\widehat{b_j-b_j+1},\cdots ,b_j-b_m+1}} \right|\, z\) = \\ \\
    \sum\limits_{i=1}^m (-1)^{m-1}   
    \prod\limits_{k=1}^m \displaystyle{\frac{\Gamma(b_k - c_i)\, \Gamma(b_j - b_k + 1)}
    {\Gamma(c_k - c_i + 1)\, \Gamma(b_j - c_k)}}\, 
    \displaystyle{\frac{e^{\pi \sqrt{-1}\, (c_i - b_j)}\, \mu_{i\, trig}^2}{\sin \pi (b_j - c_i)}}\, \times \\ \\
    z^{c_i}\! \phantom{|}_mF_{m-1}\(\left. \displaystyle{{b_1 - c_i, \cdots , b_m - c_i}
    \atop{c_1 - c_i + 1,\cdots ,\widehat{c_i - c_i + 1},\cdots ,c_m - c_i + 1}} \right|\, 
    \displaystyle{\frac1z}\).
  \end{array}
\end{equation}

\noindent Similarly, for the mHGS
\begin{equation} \label{equation:trig_continuation}
  \(T_0\)_j = \sum\limits_{i=1}^m   
  \displaystyle{\frac{\prod\limits_{k=1}^m \Gamma(b_k - c_i)}
  {\prod\limits_{k=1 \atop{k \ne i}}^m \Gamma(c_k - c_i)}}\,
  \displaystyle{\frac{\prod\limits_{k=1\atop{k \ne j}}^m \Gamma(b_j-b_k)}
  {\prod\limits_{k=1}^m \Gamma(b_j-c_k)}}\,
  \displaystyle{\frac{e^{\pi \sqrt{-1}\, (c_i-b_j)}\, \mu_{i\, trig}^2}{\sin \pi (b_j-c_i)}}\,  \(T_{\infty}\)_i. 
\end{equation}

\noindent This formula follows from Theorem \ref{thm:hg:T_0,Tinf} combined with 
(\ref{equation:analytic_cont_of_mFm-1_to_infty}). Note that all the multiples in 
(\ref{equation:g.h.g.eq.:an.cont.}) and in (\ref{equation:trig_continuation}) either bear 
the index $i$ or $j$, except for the term
$$
  \displaystyle{\frac{1}{\sin \pi (b_j - c_i)}}.
$$
So the rest of them are nothing but normalizing coefficients. Lemma \ref{lemma:p<->q} 
is true for any $a_2$, but for now $a_2 = 0$, so compairing the formula of the lemma to 
(\ref{equation:g.h.g.eq.:an.cont.}) and to (\ref{equation:trig_continuation}) 
tells us how the $p_i$ and $q_i$ are related to the hypergeometric functions of the GHGE and 
the mHGS respectively. Also, comparing (\ref{equation:trig_continuation}) to 
(\ref{equation:g.h.g.eq.:an.cont.}) is the way of going back and forth between 
the mHGS and the GHGE promised in Section \ref{sec:intro}. Finally, let us mention here that 
the explicit formulae for the local solutions (\ref{equation:T0}) and (\ref{equation:Tinfty}) 
of the mHGS similar to (\ref{equation:GHGE:sols:0}) and (\ref{equation:GHGE:sols:infty}) 
for the GHGE are given in Theorem \ref{thm:hg:T_0,Tinf}. \\
  
	Let $D_r$ ({\it r} being the first letter of the word ``rational'') be the following 
$m \times m$ matrix:
\begin{equation}
  \label{equation:D_r}
  \(D_r\)_{ij} = \displaystyle{\frac{1}{a_2 + b_i - c_j + \tau}}.
\end{equation}

\noindent The explicit formula for the inverse of this matrix 
\begin{equation} \label{equation:D_r^-1}
  \(D_r\)_{ij}^{-1} = \displaystyle{\frac{\mu_i^2(\tau)\, \nu_j^2(\tau)}
  {a_2 + b_j - c_i + \tau}}
\end{equation}

\noindent is helpful in various applications. In particular, the proof of 
Theorem \ref{thm:product} is based on it. The factorization of the determinant
\begin{equation} \label{equation:cauchy:det:rat}
  det\(\displaystyle{\frac{1}{b_i - c_j}}\) = 
  \displaystyle{\frac{\prod\limits_{1 \le i < j \le m} (b_i - b_j)(c_j - c_i)}
  {\prod\limits_{i, j = 1}^m (b_i - c_j)}}
\end{equation}

\noindent is called {\it Cauchy identity} (although known to L'Hospital and probably 
earlier). It was pointed out to the author by A.~Borodin that combining Cauchy identity 
with Cramer's rule proves (\ref{equation:D_r^-1}) in the case $a_2 = \tau = 0$. From here, 
the general case is obtained by renaming the variables. We shall call (\ref{equation:D_r^-1}) 
the {\it rational Cauchy identity} because of its equivalence, in view of Cramer's rule, to 
(\ref{equation:cauchy:det:rat}). The trigonometric analogue of 
(\ref{equation:D_r^-1}) for the matrix 
\begin{equation} \label{equation:D_trig}
  \(D_{trig}\)_{ij} = \displaystyle{\frac{1}{\sin \pi (a_2 + b_i - c_j)}}
\end{equation}

\noindent is given by a similar formula:
\begin{equation} \label{equation:D_trig^-1}
  \(D_{trig}\)_{ij}^{-1} = \displaystyle{\frac{\mu_{i\, trig}^2\, \nu_{j\, trig}^2}
  {\sin \pi (a_2 + b_j - c_i)}}.
\end{equation} 

\noindent We shall call (\ref{equation:D_trig^-1}) the {\it trigonometric Cauchy identity}.  
It is easy to obtain from the rational Cauchy identity. Start from the case $a_2 = \tau = 0$. 
Replace $b_i$ and $c_j$ by $e^{2 \pi \sqrt{-1}\, b_i}$ and $e^{2 \pi \sqrt{-1}\, c_j}$ respectively. 
Then use the identity 
$$
  e^{2 \pi \sqrt{-1}\, b_i} - e^{2 \pi \sqrt{-1}\, c_j}
   = 2\, \sqrt{-1}\, e^{\pi \sqrt{-1}\, (b_i + c_j)}\, \sin \pi (b_i - c_j),
$$ 
simplify, and finally rename the variables one last time. \\

	Let us call ${\cal M}$ and ${\cal N}$ the diagonal $m \times m$ matrices with $\mu_i$ 
and $\nu_i$ on the diagonal (we stay at the formal level and not specify the branches of 
the square roots). Then, in both the rational and the trigonometric case, we can rewrite 
the Cauchy identities as
\begin{equation} \label{equation:comp_symm}
  \({\cal N}\, D\, {\cal M}\)^t = \({\cal N}\, D\, {\cal M}\)^{-1}.
\end{equation}

\noindent The rational version of this formula lies at the core of proving the existence of 
the complex symmetric product (\ref{equation:(*,*)_tau12}) of Theorem \ref{thm:product}. 
The trigonometric Cauchy identity is used to give a proof to Theorem \ref{thm:(p,p),(q,q)} 
different from those of Beukers and Heckman and of Haraoka. \\

	Looking at (\ref{equation:D_r^-1}) and (\ref{equation:D_trig^-1}), 
one has an itch to guess the {\it elliptic Cauchy identity} for the matrix
\begin{equation} \label{equation:D_ell}
  \(D_{ell}\)_{ij} = \displaystyle{\frac{1}{sn(a_2 + b_i - c_j)}}.
\end{equation}

\noindent The problem of inverting (\ref{equation:D_ell}) turns up quite often in the theory of 
elliptic functions and, according to a specialist in the field (E.~Rains), has not been solved 
to date. We shall also need to invert (\ref{equation:D_ell}) for our own purposes. 

\begin{theorem} \label{thm:D_ell^-1}
  \begin{equation} \label{equation:D_ell^-1}
    \(D_{ell}\)_{ij}^{-1} = \displaystyle{\frac{\mu_{i\, ell}^2\, \nu_{j\, ell}^2}
    {sn(a_2 + b_j - c_i)}}
  \end{equation}
\end{theorem}

	Note that the trigonometric and elliptic versions of the determinantal identity 
(\ref{equation:cauchy:det:rat}) are easy to obtain from (\ref{equation:D_trig^-1}) and 
(\ref{equation:D_ell^-1}) by means of Cramer's rule.

\subsection{Fermionic fields}% %%%%%%%%%%%%%%%%%%%%%%%%%%%%%%%%%

	Let $\omega_1$ and $\omega_2$ be a pair of linearly independent complex 
numbers and let $\Omega = \left\{n_1 \omega_1 + n_2 \omega_2 : n_1,n_2 \in \Z \right\}$. 
From this point on, let $\tau_1 = n_1 \omega_1$ and $\tau_2 = n_2 \omega_2$. \\ 

	Consider the following vectors in $H_{\tau_1, \tau_2} \oplus H_{\tau_2, \tau_1}$:
\begin{equation}
\label{equation:f}
\begin{array}{l}
\(f_0\)_i(\tau_1, \tau_2) = \displaystyle{\frac{1}{\nu_i^2(\tau_1 + \tau_2)}}\, 
\displaystyle{\frac{v_i(\tau_1) \oplus \sqrt{-1}\, v_i(\tau_2)}{\sqrt{2}}}, \\ \\
\(f_0\)_i^{\dagger}(\tau_1, \tau_2) = \displaystyle{\frac{1}{\nu_i^2(\tau_1 + \tau_2)}}\, 
\displaystyle{\frac{v_i(\tau_1) \oplus -\sqrt{-1}\, v_i(\tau_2)}{\sqrt{2}}}, \\ \\
\(f_{\infty}\)_i(\tau_1, \tau_2) = \displaystyle{\frac{1}{\mu_i^2(\tau_1 + \tau_2)}}\, 
\displaystyle{\frac{w_i(\tau_2) \oplus \sqrt{-1}\, w_i(\tau_1)}{\sqrt{2}}}, \\ \\ 
\(f_{\infty}\)_i^{\dagger}(\tau_1, \tau_2) = \displaystyle{\frac{1}{\mu_i^2(\tau_1 + \tau_2)}}\, 
\displaystyle{\frac{w_i(\tau_2) \oplus -\sqrt{-1}\, w_i(\tau_1)}{\sqrt{2}}}.
\end{array}
\end{equation}

\noindent Let us introduce a modified symmetric scalar product on
$H_{\tau_1, \tau_2} \oplus H_{\tau_2, \tau_1}$. For $\tau_1 + \tau_2 \ne 0$, 
\begin{equation} \label{equation:scalprod:tau_1tau_2}
  \begin{array}{l}
    \(g_1 \oplus g_2, h_1 \oplus h_2\) = (-1)^{n_1} 
    \(\(g_1, h_1\)_{\tau_1, \tau_2} + \(g_2, h_2\)_{\tau_2, \tau_1} \right. - \\ \\ 
    \left. \displaystyle{\frac{1}{\tau_1 + \tau_2}} 
    \(\(u, g_1\)_{\tau_1, \tau_2} \(u, h_1\)_{\tau_1, \tau_2} + 
    \(u, g_2\)_{\tau_2, \tau_1} \(u, h_2\)_{\tau_2, \tau_1}\)\). 
  \end{array}
\end{equation}

\noindent For $\tau_1 + \tau_2 = 0$, let us define
\begin{equation}
\label{equation:scalprod:tau_1tau_2:0}
\(g_1 \oplus g_2, h_1 \oplus h_2\) = \(g_1, h_1\)_{\tau_1, \tau_2} + 
\(g_2, h_2\)_{\tau_2, \tau_1}.
\end{equation}

\noindent Here is the multiplication table for (\ref{equation:scalprod:tau_1tau_2}):
\begin{equation}
\label{equation:mult_table}
\begin{array}{|c|}
\hline \\
\(\(f_0\)_i(\tau_1, \tau_2), \(f_0\)_j(\tau_1, \tau_2)\) = 
\(\(f_0\)_i^{\dagger}(\tau_1, \tau_2), 
\(f_0\)_j^{\dagger}(\tau_1, \tau_2)\) = 0 \\ \\
\hline \\
\(\(f_0\)_i(\tau_1, \tau_2), \(f_0\)_j^{\dagger}(\tau_1,\tau_2)\) = 
(-1)^{n_1} \(\displaystyle{\frac{\delta_{ij}}{\nu_i^2(\tau_1 + \tau_2)}} - 
\displaystyle{\frac{1}{\tau_1 + \tau_2}}\) \\ \\ \hline \\
\(\(f_{\infty}\)_i(\tau_1, \tau_2), \(f_{\infty}\)_j(\tau_1, \tau_2)\) =  
\(\(f_{\infty}\)_i^{\dagger}(\tau_1, \tau_2) 
\(f_{\infty}\)_j^{\dagger}(\tau_1, \tau_2)\) = 0 \\ \\ \hline \\
\(\(f_{\infty}\)_i(\tau_1, \tau_2), \(f_{\infty}\)_j^{\dagger}(\tau_1, \tau_2)\) = 
(-1)^{n_1} \(\displaystyle{\frac{\delta_{ij}}{\mu_i^2(\tau_1 + \tau_2)}} - 
\displaystyle{\frac{1}{\tau_1 + \tau_2}}\) \\ \\ \hline \\
\(\(f_0\)_i(\tau_1, \tau_2), \(f_{\infty}\)_j^{\dagger}(\tau_1, \tau_2)\) = 
\(\(f_0\)_i^{\dagger}(\tau_1, \tau_2), \(f_{\infty}\)_j(\tau_1, \tau_2)\) = \\ \\
(-1)^{n_1} \(\displaystyle{\frac{1}{a_2 + b_i - c_j + \tau_1 + \tau_2}} - 
\displaystyle{\frac{1}{\tau_1 + \tau_2}}\) \\ \\
\hline
\end{array}
\end{equation}
\medskip

	Let $H$ be a complex vector space of even dimension endowed with 
a non-degenerate symmetric scalar product $(*,*)$. A subspace $I \subset H$ 
is called {\it isotropic}, if $(v,v) = 0$ for any $v \in I$. Let $H = I \oplus I^{\dagger}$
be a decomposition of $H$ into a direct sum of maximal isotropic subspaces. Let us 
choose bases $v_j$ and $v_i^{\dagger}$ of $I$ and $I^{\dagger}$ respectively. 
Then $(v_i,v_j) = (v_i^{\dagger},v_j^{\dagger}) = 0$. It is customary to take the dual 
bases for $I$ and $I^{\dagger}$ so that $(v_i, v_j^{\dagger}) = \delta_{ij}$, but we shall 
not do so in this paper. The vectors $v_i$ and $v_i^{\dagger}$ are called 
{\it annihilation operators} and {\it creation operators} respectively. Both 
the annihilation and creation operators also bear the common name 
of {\it fermions}. \\

	A {\it Clifford algebra} $CA$ is the associative algebra generated by 
the vectors of $H$ with relations
$$
f h +h f = (f,h).
$$

\noindent For a Clifford algebra $CA$, let us call $Ann$ and $Cr$ the spaces of annihilation
and creation operators respectively. Then the left and right $CA$-modules 
$Fock = CA/CA\, Ann$ and $Fock^{\dagger} = Cr\, CA \backslash CA$ are called
the {\it Fock space} and the {\it dual Fock space} respectively.
The generators $1\, mod\, CA\, Ann$ and $1\, mod\, Cr\, CA$ are denoted by 
$|0\>$ and $\<0|$ and called the {\it vacuum vector} and the {\it dual vacuum vector}. 
The spaces $Fock$ and $Fock^{\dagger}$ are dual via the bilinear form 
$\(\<0|f,h|0\>\) \mapsto \<fh\>$ where
\begin{equation} \label{equation:Wick's rule}
\begin{array}{l}
  \<1\> = 1; \\
   \\
  \<fh\> = (f,h),~~~\mbox{if}~f,h \in W; \\
   \\
  \<h_1 \cdots h_r\> = 
    \left\{ \begin{array}{lr} 
      0, & \mbox{if}~r~\mbox{is odd}, \\
       & \\ 
      \sum\limits_{\sigma} \mbox{sign}(\sigma) \<h_{\sigma(1)} h_{\sigma(2)}\>
      \cdots \<h_{\sigma(r-1)} h_{\sigma(r)}\>, & \mbox{if}~r~\mbox{is even}.
    \end{array}\right.
\end{array}
\end{equation}
The sum runs over all the permutations $\sigma$ satisfying 
$\sigma(1) < \sigma(2), \cdots \sigma(r-1) < \sigma(r)$ and 
$\sigma(1) < \sigma(3) < \cdots < \sigma(r-1)$, in other words, 
over all ways of grouping the $h_i$ into pairs. The equation
(\ref{equation:Wick's rule}) is called {\it Wick's rule}. The number
$\<h_1 \cdots h_r\>$ is called the {\it vacuum expectation value}. \\

	Let $H = \bigoplus_{\tau_1 + \tau_2 \in \Omega} H_{\tau_1, \tau_2} \oplus 
H_{\tau_2, \tau_1}$ equipped with the scalar product (\ref{equation:scalprod:tau_1tau_2}). 
Let $CA$ be the Clifford algebra generated by $H$. Consider the following fermionic 
fields:	
\begin{equation} \label{equation:generating_functions}
  \begin{array}{l}
    \(F_0\)_i (z_1, z_2) = \sum\limits_{\tau_1 + \tau_2 \in \Omega} 
    \(f_0\)_i(\tau_1, \tau_2)\, z_1^{n_1}\, z_2^{n_2}, \\ \\ 
    \(F_0\)_i^{\dagger} (z_1, z_2) = \sum\limits_{\tau_1 + \tau_2 \in \Omega} 
    \(f_0\)_i^{\dagger}(\tau_1, \tau_2)\, z_1^{n_1}\, z_2^{n_2}, \\ \\
    \(F_{\infty}\)_i (z_1, z_2) = \sum\limits_{\tau_1 + \tau_2 \in \Omega} 
    \(f_{\infty}\)_i (\tau_1, \tau_2)\, z_1^{-n_1}\, z_2^{-n_2}, \\ \\  
    \(F_{\infty}\)_i^{\dagger} (z_1, z_2) = \sum\limits_{\tau_1 + \tau_2 \in \Omega} 
    \(f_{\infty}\)_i^{\dagger}(\tau_1, \tau_2)\, z_1^{-n_1}\, z_2^{-n_2}.
  \end{array}
\end{equation}

\begin{theorem}
  \label{thm:<F_inftyF_0>:ell}
  $$
    \left\< \(F_{\infty}\)_j^{\dagger}\, \(F_0\)_i \right\> = 
    \displaystyle{\frac{1}{sn(a_2 + b_i - c_j)}} =
    \left\< \(F_{\infty}\)_j\, \(F_0\)_i^{\dagger} \right\>,
  $$
  \smallskip
  
  \noindent where $<*>$ is the vacuum expectation value.
\end{theorem}

\proof It immediately follows from the last formula of the multiplication table 
(\ref{equation:mult_table}) that
\begin{equation}
\label{equation:<F_inftyF_0>:elliptic}
\left\< \(F_{\infty}\)_j^{\dagger}\, \(F_0\)_i\right\> = 
\displaystyle{\frac{1}{a_2 + b_i - c_j}} +
\sum\limits_{\tau_1 + \tau_2 \in \Omega \setminus\{0\}} (-1)^{n_1} \(
\displaystyle{\frac{1}{a_2 + b_i - c_j + \tau_1 + \tau_2}} - \displaystyle{\frac{1}{\tau_1 + \tau_2}}\).
\end{equation}

\noindent The following decomposition
\begin{equation} \label{equation:sn<->sin}
  \sqrt{{\cal P}(u) - e_3} = \displaystyle{\frac{1}{\omega_1}} \sum\limits_{n_2 \in \Z}\, 
  \displaystyle{\frac{\pi}{\sin \pi \(\frac{u}{\omega_1} + \frac{n_2 \omega_2}{\omega_1}\)}}
\end{equation}

\noindent is proven in different notations (see below) at the end of Chapter 2 of \cite{HC}. 
Here ${\cal P}(z)$ is the ${\cal P}$-function of Weierstrass and $e_3 = {\cal P}\(\omega_2/2\)$. 
Expanding (\ref{equation:sn<->sin}) further by means of the famous identity 
\begin{equation} \label{equation:sin}
  \displaystyle{\frac{\pi}{\sin \pi z}} = \displaystyle{\frac{1}{z}} +
  \sum\limits_{n_1 = 1}^{\infty} (-1)^{n_1} \(
  \displaystyle{\frac{1}{z + n_1}} + \displaystyle{\frac{1}{z - n_1}}\),
\end{equation}

\noindent produces the following formula:
\begin{equation} \label{equation:1/sn:sum}
  \begin{array}{l}
     \displaystyle{\sqrt{{\cal P}(a_2 + b_i - c_j) - e_3}} = \\ \\ 
     \displaystyle{\frac{1}{a_2 + b_i - c_j}} + \sum\limits_{\tau_1 + \tau_2 \in \Omega \setminus\{0\}} 
    (-1)^{n_1} \(\displaystyle{\frac{1}{a_2 + b_i - c_j + \tau_1 + \tau_2}} -
    \displaystyle{\frac{1}{\tau_1 + \tau_2}}\).
  \end{array}
\end{equation}

\noindent The last term appears in (\ref{equation:1/sn:sum}) to make the sum 
converge uniformly. Similarly, one can rewrite (\ref{equation:sin}) as
$$
  \displaystyle{\frac{\pi}{\sin \pi z}} = \displaystyle{\frac{1}{z}} +
  \sum\limits_{n_1 \in Z \setminus\{0\}}^{\infty} (-1)^{n_1} \(
  \displaystyle{\frac{1}{z + n_1}} - \displaystyle{\frac{1}{n_1}}\).
$$

	The elliptic sine of Jacobi $sn(z)$ is sometimes defined as
\begin{equation}
\label{equation:sn}
\displaystyle{\frac{1}{sn(z)}} = \sqrt{{\cal P}(z) - e_3}.
\end{equation}

\noindent It is customary in the theory of Jacobi elliptic functions to use notations
different from those of Weierstrass:
$$
\begin{array}{llll}
\omega = \displaystyle{\frac{\omega_1}{2}}, & 
\omega' = \displaystyle{\frac{\omega_2}{2}}, & 
\tau = \displaystyle{\frac{\omega'}{\omega}} = 
\displaystyle{\frac{\omega_2}{\omega_1}}, & 
h = e^{\pi \sqrt{-1}\, \tau}. 
\end{array}
$$ 

\noindent The generators $4\omega$ and $2\omega'$ of the periods lattice are
chosen so that $\mbox{Im}\, \tau > 0$. \\

	It is standard to add the following normalizing condition to 
the definition (\ref{equation:sn}) of $sn(z)$:
\begin{equation} \label{equation:normSn}
  \omega = \displaystyle{\frac{\pi}{2}} \(\sum\limits_{n = -\infty}^{\infty} h^{n^2}\)^2,
\end{equation}

\noindent and to treat $sn(z) = sn(z; \tau)$ rather than $sn(z) = sn(z; \omega, \omega')$;
see \cite{HC} or \cite{WW} for more information on elliptic functions. We shall use
the definition (\ref{equation:sn}) of $sn(z)$ without the condition (\ref{equation:normSn}). 
The only reason we switch notations form $\cal{P}$ to $sn$ is that the latter takes less
space. This remark concludes the proof. $\Box$ \\

	Let us call $\H_0$,  $\H_0^{\dagger}$, $\H_{\infty}$, and $\H_{\infty}^{\dagger}$ 
the $m$-dimensional spaces spanned by $\(F_0\)_i$, $\(F_0\)_i^{\dagger}$, 
$\(F_{\infty}\)_i$, and $\(F_{\infty}\)_i^{\dagger}$ respectively. We use the pairing of 
Theorem \ref{thm:<F_inftyF_0>:ell} to identify $\H_0$ with $\H_{\infty}^{\dagger}$ and
$\H_0^{\dagger}$ with $\H_{\infty}$. The resulting vector spaces are $\H$ and $\H'$ 
of Theorem \ref{thm:(F,F):ell}. 

\subsection{Trigonometrization}%%%%%%%%%%%%%%%%%%%%%%%%%%%%%%%

Let us restrict ourselves to the case $\omega_1 = 1$ and let us consider what happens 
when $\omega_2 \to \infty$. The limits of $\(f_0\)_i(\tau_1, \tau_2)$, 
$\(f_0\)_i^{\dagger}(\tau_1, \tau_2)$, $\(f_{\infty}\)_i(\tau_1, \tau_2)$, and 
$\(f_{\infty}\)_i^{\dagger}(\tau_1, \tau_2)$ are all zero unless $\tau_2 = 0$. So in the limit 
we get 
\begin{equation} \label{equation:generating_functions:trig}
  \begin{array}{l}
    \lim\limits_{\omega_2 \to \infty} \(F_0\)_i(z_1, z_2) = \(F_0\)_{i\, trig}(z) = 
    \sum\limits_{n \in \Z} \(f_0\)_i(n,0)\, z^n, \\ \\
    \lim\limits_{\omega_2 \to \infty} \(F_0\)_i^{\dagger}(z_1, z_2) = \(F_0\)_{i\, trig}^{\dagger}(z) = 
    \sum\limits_{n \in \Z} \(f_0\)_i^{\dagger}(n, 0)\, z^n, \\ \\
    \lim\limits_{\omega_2 \to \infty} \(F_{\infty}\)_i(z_1, z_2) = \(F_{\infty}\)_{i\, trig}(z) = 
    \sum\limits_{n \in \Z} \(f_{\infty}\)_i(n,0)\, z^{-n}, \\ \\
    \lim\limits_{\omega_2 \to \infty} \(F_{\infty}\)_i^{\dagger}(z_1, z_2) = 
    \(F_{\infty}\)_{i\, trig}^{\dagger}(z) = \sum\limits_{n \in \Z} \(f_{\infty}\)_i^{\dagger}(n, 0)\, z^{-n}.
  \end{array}
\end{equation}

\noindent This time the vacuum expectation values are
\begin{equation} \label{equation:<F_inftyF_0>:trig}
  \begin{array}{l}
    \left\< \(F_{\infty}\)_{j\, trig}^{\dagger}\, \(F_0\)_{i\, trig} \right\> = 
    \left\< \(F_{\infty}\)_{j\, trig}\, \(F_0\)_{i\, trig}^{\dagger} \right\> = \\ \\
    \displaystyle{\frac{1}{a_2 + b_i - c_j}} + \sum\limits_{n \in \Z \setminus \{0\}} 
    (-1)^n \(\displaystyle{\frac{1}{a_2 + b_i - c_j + n}} - \displaystyle{\frac1n}\) = \\ \\
    \displaystyle{\frac{\pi}{\sin \pi (a_2 + b_i - c_j)}} = 
    \lim\limits_{\omega_2 \to \infty} \displaystyle{\frac{1}{sn(a_2 + b_i - c_j)}}.
  \end{array}
\end{equation}

	Let us call $H_{trig} = \bigoplus_{n \in \Z} H_{n0} \oplus H_{0n}$ 
and $CA_{trig}$ the Clifford algebra generated by $H_{trig}$. We shall call 
$\H_{0\, trig}$,  $\H_{0\, trig}^{\dagger}$, $\H_{\infty\, trig}$, and $\H_{\infty\, trig}^{\dagger}$ 
the $m$-dimensional spaces spanned by $\(F_0\)_{i\, trig}$, $\(F_0\)_{i\, trig}^{\dagger}$, 
$\(F_{\infty}\)_{i\, trig}$, and $\(F_{\infty}\)_{i\, trig}^{\dagger}$ respectively. We use 
the pairing (\ref{equation:<F_inftyF_0>:trig}) to identify $\H_{0\, trig}$ with 
$\H_{\infty\, trig}^{\dagger}$. We call the resulting vector space $\H_{trig}$. Similarly, 
we identify $\H_{0\, trig}^{\dagger}$ with $\H_{\infty\, trig}$ and call the resulting space 
$\H'_{trig}$. The following theorem is a trigonometric version of Theorems \ref{thm:(F,F):ell} 
and \ref{thm:hk} combined:
\begin{theorem} \label{thm:herm@F}
  Let the local exponents of the mHGS $b_1, \dots , b_m$ and $-c_1, \dots ,-c_m$ be generic 
  real numbers. Then there exist unique up to a constant multiple hermitian and complex 
  symmetric scalar products $(*,*)_{trig}$ on $\H_{trig}$ such that the bases 
  $\(F_0\)_{i\, trig}$ and $\(F_{\infty}\)_{i\, trig}^{\dagger}$ are simultaneously 
  orthogonal with respect to them, given by 
  \begin{equation} \label{equation:trig:length}
    \(\(F_0\)_{i\, trig}, \(F_0\)_{j\, trig}\)_{trig} = \displaystyle{\frac{\delta_{ij}}{\nu_{i\, trig}^2}} 
    \mbox{~~and~~} 
    \(\(F_{\infty}\)_{i\, trig}^{\dagger}, \(F_{\infty}\)_{j\, trig}^{\dagger}\)_{trig} = 
    \displaystyle{\frac{\delta_{ij}}{\mu_{i\, trig}^2}}.
  \end{equation}

  \noindent A similar statement holds for $\H'_{trig}$. \\

    The quaternions act on $\H_{trig} \oplus \H'_{trig}$ according to (\ref{equation:quaternions:F}).
  If the forms are sign-definite, then the quaternionic action is hyperk\"ahler.
\end{theorem}
\medskip

	The map 
$$
  p_i \mapsto \nu_{i\, trig}^2 \(F_0\)_{i\, trig}, ~~~ 
  q_i \mapsto \mu_{i\, trig}^2 \(F_{\infty}\)_{i\, trig}^{\dagger}
$$
establishes hermitian isometries between $\H_{trig}$ and the spaces of solution of the GHGE 
and the mHGS, see Theorem \ref{thm:(p,p),(q,q)} and Lemma \ref{lemma:p<->q}. Let us repeat 
here that the isometry is a formal linear-algebraic construction. It allows to introduce the 
quaternionic action on the direct sum of the spaces of solutions and their duals, but the action 
of the monodromy group on the fermionic fields doesn't make geometric sense. The situation 
should be improved by constructing a Fuchsian system on the torus 
$T(\omega_1, \omega_2)$ such that as $\omega_2 \to \infty$, the system becomes the mHGS. 

%%%%%%%%%%%%%%%%%%%%%%%%%%%%%%%%%%%%%%%%%%%%%%%%
\section{Proofs}%%%%%%%%%%%%%%%%%%%%%%%%%%%%%%%%%%%%%%%%
In this section, we prove the above results. The section also contains some more new results, 
which we use as tools, but which we believe are important for their own sake. These are 
Theorem \ref{thm:hg:T_0,Tinf}, the extended version of Theorem \ref{thm:product} on page 
\pageref{thm:product:ext} and an explicit construction of the monodromy operators 
of the mHGS as well as their eigenvectors in a ``good'' basis in Subsection \ref{subsec:BH}.  

\subsection{Rigidity and irreducibility}%%%%%%%%%%%%%%%%%%%%%%%%%%%%%%

Let $R_1, \dots ,R_k$ ($M_1, \dots ,M_k$) be a $k$-tuple of linear operators on $\C^m$ 
such that $R_1 + \dots +R_k = 0$ \\ ($M_1 \cdot \dots \cdot M_k = Id$, in which case 
the $M_i$ are called a {\it local system}). The $k$-tuple is called {\it irreducible}, if the operators 
do not simultaneously preserve a proper subspace of $\C^m$. Let $R'_1, \dots ,R'_k$ 
($M'_1, \dots ,M'_k$) be any other $k$-tuple of linear operators on $\C^m$ such that 
$R'_1 + \dots +R'_k = 0$ \\ ($M'_1 \cdot \dots \cdot M'_k = Id$) and that the operators 
$R_i$ and $R'_i$ ($M_i$ and $M'_i$) are conjugate to each other for $i = 1, \dots ,k$. 
The tuple is called {\it rigid} \label{def:rigid}, if all the operators $R_i$ and $R'_i$ 
($M_i$ and $M'_i$) are simultaneously conjugate to each other. The additive and 
multiplicative hypergeometric triples are rigid and irreducible (see \cite{BH}, \cite{Ko1}, 
and \cite{S} for proofs). \\

	Consider the following matrices:
\begin{equation} \label{equation:A,B,C}
  \begin{array}{l}
    A_{ij} = \left\{
    \begin{array}{lll}
      c_i - b_i         & ,if & i=j \\
      c_i - b_i -a_2 & ,if & i \ne j 
    \end{array} \right. , \\ \\

    B_{ij} = \left\{
    \begin{array}{lll}
      a_2 + b_i - c_i & ,if & i<j \\
      b_i                   & ,if & i=j \\
      0                      & ,if & i>j
    \end{array} \right. , \\ \\
    
    C_{ij} = \left\{
    \begin{array}{lll}
      0                      & ,if & i<j \\
      -c_i                  & ,if & i=j \\
      a_2 + b_i - c_j & ,if & i>j
    \end{array} \right. .
  \end{array}
\end{equation}

\noindent Quite obviously, they add up to zero and have the same eigenvalues 
(with multiplicities) as the operators $A$, $B$, and $C$ we have started with. It
is proven in \cite{G1} that the following vectors are the eigenvectors of $B$ and $C$
corresponding to the eigenvalues $b_i$ and $-c_i$ respectively:
\begin{equation} \label{equation:v} 
  v^j_i = \left\{
  \begin{array}{cl}
    (b_j - c_j + a_2)\, \displaystyle{\frac{\prod\limits_{k = j + 1}^m (b_i - c_k + a_2)}
    {\prod\limits_{ k = j \atop{k \ne i}}^m (b_i - b_k)}}, & if~~~ i \ge j, \\
     & \\
    0, & if~~~ i < j; 
  \end{array} \right.
\end{equation}
\medskip

\begin{equation} \label{equation:w}
  w^j_i = \left\{
  \begin{array}{cl} 
    0, & if~~~ i > j; \\
     & \\
    (b_j - c_j + a_2)\, \displaystyle{\frac{\prod\limits_{k=1}^{j - 1} (b_k - c_i + a_2)}
    {\prod\limits_{k=1 \atop{k \ne i}}^j (c_k - c_i)}}, & if~~~ i \le j. 
  \end{array} \right.
\end{equation}

\noindent Here and in the sequel, all empty products are understood to be equal to $1$. \\ 

	Let ${\cal F}_0$ and ${\cal F}_{\infty}$ be the flags 
$Span(v_1) \subset Span(v_1, v_2) \subset \dots \subset \C^m$ and 
$Span(w_m) \subset Span(w_m, w_{m-1}) \\ \subset \dots \subset \C^m$. The flags are 
opposite to each other as 
${\cal F}_0 = Span(e_1) \subset Span(e_1, e_2) \subset \dots \subset \C^m$, whereas 
${\cal F}_{\infty} = Span(e_m) \subset Span(e_m, e_{m-1}) \subset \dots \subset \C^m$. 
Thus, the flags ${\cal F}_0$ and ${\cal F}_{\infty}$ are in general position with respect to 
each other: if we take a subspace from ${\cal F}_0$ and from ${\cal F}_{\infty}$, 
the dimension of their intersection will be the lowest possible. \\

	Let $u = (b_1 - c_1 + a_2, \dots ,b_m - c_m + a_2)$. 
It is not hard to see that $A u = a_1 u$. It is proven in \cite{G1} that 
$$ \sum\limits_{i = 1}^m v_i = \sum\limits_{i = 1}^m w_i = u. $$ Thus the flag 
${\cal F}_1 = Span(u) \subset \C^m$ is in general position with respect to the flags 
${\cal F}_0$ and ${\cal F}_{\infty}$. The three {\it spectral flags} being in general position 
with respect to each other is the combinatorial part of the condition that the operators 
$A$, $B$, and $C$ are generic (see also \cite{Mag1}). Due to rigidity, any additive 
hypergeometric triple of matrices with the same eigenvalues is conjugate to the triple above. 
This fact allows us to make most of the proofs computational. 

\subsection{Main technical tools}%%%%%%%%%%%%%%%%%%%%%%%%%%%%%%
For $\tau \in \C^*$, let
\begin{equation} \label{equation:Xi,Theta}
  \begin{array}{lll}
    \xi_i^2(\tau) = \displaystyle{\frac{\prod\limits_{k=1}^m 
    (b_i - b_k + \tau)}{\prod\limits_{k=1\atop{k \ne i}}^m (b_i - b_k)}}, & &
    \theta_i^2(\tau) = \displaystyle{\frac{\prod\limits_{k=1}^m (c_k - c_i + \tau)}
    {\prod\limits_{k=1\atop{k \ne i}}^m (c_k - c_j)}}, \\
     & \phantom{Hi there!} & \\
    \Xi^2(\tau) = diag\(\xi_i^2(\tau)\)_{i=1, \cdots ,m}, & &
    \Theta^2(\tau) = diag\(\theta_i^2(\tau)\)_{i=1, \cdots ,m}.
  \end{array}
\end{equation}

\noindent For $\tau \in \C^*$, let $EX(\tau)$ and $EY(\tau)$ be the following 
$m \times m$ matrices:
\begin{equation} \label{equation:EX,EY}
  EX_{ij}(\tau) = \displaystyle{\frac{\xi_j^2(\tau)}{b_j - b_i + \tau}},~~~
  EY_{ij}(\tau) = \displaystyle{\frac{\theta_j^2(\tau)}{c_i - c_j + \tau}}.
\end{equation}

\noindent As $\lim\limits_{\tau \to 0} EX(\tau) = \lim\limits_{\tau \to 0} EY(\tau) = Id$, 
we naturally set $EX(0) = EY(0) = Id$. As we shall show in Lemma \ref{lemma:expX,Y}, 
in fact  
$$
EX(\tau) = e^{X \tau} \mbox{~~and~~} EY(\tau) = e^{Y \tau},
$$
where $X$ and $Y$ are the matrices (\ref{equation:X}) and (\ref{equation:Y}) respectively. 
\newpage

\begin{lemma} \label{lemma:EX,EY} $\phantom{tra-la-la}$ \\

\noindent $\bullet$~~~ For any $\tau_1, \tau_2 \in \C$, 
$EX(\tau_1)\, EX(\tau_2) = EX(\tau_1 + \tau_2)$ and 
$EY(\tau_1)\, EY(\tau_2) = EY(\tau_1 + \tau_2)$. \\

\noindent $\bullet$~~~  For $\tau \ne 0$, the Jordan normal form of $EX(\tau)$ 
and $EY(\tau)$ is a single block with the eigenvalue $1$. \\

\noindent $\bullet$~~~ Recall that $e=(1, \cdots ,1)$. $EX(\tau)\, e = 
EY(\tau)\, e = e$.
\end{lemma}

\proof All the proofs for $EX$ and $EY$ are the same, so we shall only prove 
the statements of the lemma for $EX$. \\

\noindent $\bullet$~~~ To prove the first statement of the lemma, we need to show that 
$\(EX(\tau_1)\, EX(\tau_2)\)_{ij} = EX(\tau_1 + \tau_2)_{ij}$ i.e.
\begin{equation} \label{equation:lemma:EX:1}
  \sum\limits_{k=1}^m \displaystyle{\frac{\prod\limits_{l=1 \atop{l \ne i}}^m (b_k - b_l + \tau_1)}
  {\prod\limits_{l=1 \atop{l \ne k}}^m (b_k - b_l)}}\, 
  \displaystyle{\frac{\prod\limits_{l=1 \atop{l \ne k}}^m (b_j - b_l + \tau_2)}
  {\prod\limits_{l=1 \atop{l \ne j}}^m (b_j - b_l)}} =
  \displaystyle{\frac{\prod\limits_{l=1 \atop{l \ne i}}^m (b_j - b_l + \tau_1 + \tau_2)}
  {\prod\limits_{l=1 \atop{l \ne j}}^m (b_j - b_l)}}.
\end{equation}
Canceling out common multiples, we rewrite (\ref{equation:lemma:EX:1}) as
\begin{equation} \label{equation:lemma:EX:2}
  \sum\limits_{k=1}^m \displaystyle{\frac{\prod\limits_{l=1 \atop{l \ne i}}^m (b_k - b_l + \tau_1)}
  {\prod\limits_{l=1 \atop{l \ne k}}^m (b_k - b_l)}}\, 
  \displaystyle{\frac{\prod\limits_{l=1 \atop{l \ne k}}^m (b_j - b_l + \tau_2)}
  {\prod\limits_{l=1 \atop{l \ne i}}^m (b_j - b_l + \tau_1 + \tau_2)}} = 1.
\end{equation}

	Rational identities of this kind will appear often in the rest of the paper. Let us outline 
the strategy of proving them here. The left hand side $L$ of an identity to prove
will be a homogeneous rational function of $b_i$, $c_i$, $\tau_1$, and $\tau_2$. 
The right hand side $R$ of the identity will be a homogeneous polynomial in the same 
variables such that $deg(L) = deg(R)$. The degree will not exceed $1$. 
All the denominators of $L$ will be products of linear forms $f$. The power of 
every such form in every denominator will be $1$. The first step to prove such an identity 
is to show that $L$ is in fact a polynomial. For that, it is enough to show that 
$f\, L|_{f=0}=0$ for every form $f$ from any denominator of the identity. The second step is 
to check enough points to make sure that $L \equiv R$. We shall write down the proof 
of (\ref{equation:lemma:EX:2}) in full detail. We shall be sketchy with the rest of 
the proofs, if they follow the strategy outlined here. \\

Let $1 \le l_1 < l_2 \le m$. The only two summands in $(b_{l_1} - b_{l_2}) L$
which do not nullify when restricted to $b_{l_1} = b_{l_2}$ are 
$$
\displaystyle{\frac{b_{l_1} - b_{l_2}}{b_{l_1} - b_{l_2}}}\, 
\displaystyle{\frac{\prod\limits_{l=1 \atop{l \ne i}}^m (b_{l_1} - b_l + \tau_1)}
{\prod\limits_{l=1 \atop{l \ne l_1 \atop{l \ne l_2}}}^m (b_{l_1} - b_l)}}\, 
\displaystyle{\frac{\prod\limits_{l=1 \atop{l \ne l_1}}^m (b_j - b_l + \tau_2)}
{\prod\limits_{l=1 \atop{l \ne i}}^m (b_j - b_l + \tau_1 + \tau_2)}}
~\mbox{and}~
\displaystyle{\frac{b_{l_1} - b_{l_2}}{b_{l_2} - b_{l_1}}}\, 
\displaystyle{\frac{\prod\limits_{l=1 \atop{l \ne i}}^m (b_{l_2} - b_l + \tau_1)}
{\prod\limits_{l=1 \atop{l \ne l_2 \atop{l \ne l_1}}}^m (b_{l_2} - b_l)}}\, 
\displaystyle{\frac{\prod\limits_{l=1 \atop{l \ne l_2}}^m (b_j - b_l + \tau_2)}
{\prod\limits_{l=1 \atop{l \ne i}}^m (b_j - b_l + \tau_1 + \tau_2)}}.
$$

\noindent Obviously, these two add up to zero. \\

	For $l_1 = 1, \cdots ,{\hat i}, \cdots ,m$, let $f = b_j - b_{l_1} + \tau_1 + \tau_2$.
This time
$$
\begin{array}{l}
  f\, L|_{f = 0} = 
  \left. \displaystyle{\frac{\prod\limits_{l=1}^m (b_j - b_l + \tau_2)}
  {\prod\limits_{l=1 \atop{l \ne i \atop{l \ne l_1}}}^m (b_j - b_l + \tau_1 + \tau_2)}}\,
  \sum\limits_{k=1}^m \displaystyle{\frac{\prod\limits_{l=1 \atop{l \ne i}}^m (b_k - b_l + \tau_1)}
  {\prod\limits_{l=1 \atop{l \ne k}}^m (b_k - b_l)}}\, \displaystyle{\frac{1}
  {b_j - b_k + \tau_2}}\right|_{b_j - b_{l_1} + \tau_1 + \tau_2 = 0} = \\ \\
  \displaystyle{\frac{\prod\limits_{l=1}^m (b_l - b_{l_1} + \tau_1)}
  {\prod\limits_{l=1 \atop{l \ne i \atop{l \ne l_1}}}^m (b_l - b_{l_1})}}\,
  \sum\limits_{k=1}^m \displaystyle{\frac{\prod\limits_{l=1 \atop{l \ne i}}^m (b_k - b_l + \tau_1)}
  {\prod\limits_{l=1 \atop{l \ne k}}^m (b_k - b_l)}}\, \displaystyle{\frac{-1}{b_k - b_{l_1} + \tau_1}}.
\end{array}
$$

\noindent The identity 
$$
\sum\limits_{k=1}^m \displaystyle{\frac{\prod\limits_{l=1 \atop{l \ne i \atop{l \ne l_1}}}^m 
(b_k - b_l + \tau_1)}{\prod\limits_{l=1 \atop{l \ne k}}^m (b_k - b_l)}} = 0
$$
for $l_1 \ne i$ is proven in \cite{G1}, see (7.59) of the Appendix there. \\

	We now know that $L$ is a homogeneous polynomial in $b_i$ and $\tau_i$. 
The degree of this polynomial is zero, so it must be a constant. Setting $b_k = k$ for
$k = 1, \cdots, m$ and $\tau_1 = 0$, we see that all the summands
on the left hand side of (\ref{equation:lemma:EX:2}) nullify except for when $k=i$. 
The one remaining equals the one on the opposite side of the identity. \\

\noindent $\bullet$~~~ To prove the second statement of the lemma, let us introduce 
the following matrix:
\begin{equation} \label{equation:G}
  \begin{array}{r}
    G(\tau)_{ij} = \displaystyle{\frac{1}{\prod\limits_{k = 1 \atop{k \ne j}}^m (b_j - b_k)}}
    \[1 + \displaystyle{\frac{b_j - b_1}{\tau}} + 
    \displaystyle{\frac{(b_j - b_1)(b_j - b_1 - \tau)}{2!\, \tau^2}} + \cdots \right. \\
    \left. \cdots + 
    \displaystyle{\frac{(b_j - b_1) \cdot \ldots \cdot (b_j - b_1 - (m - i -1) \tau)}
    {(m - i)!\, \tau^{m - i}}}\].
  \end{array}
\end{equation}

\noindent To prove that $G(\tau)$ in not degenerate for $\tau \ne 0$, 
let us prove that  
\begin{equation} \label{equation:det(G)}
  \det\(G(\tau)\) = \displaystyle{\frac{1}{1!\, 2!\, \ldots (m-1)!\, \tau^{m(m - 1)/2}\, 
  \prod\limits_{1 \le i < j \le m} (b_j - b_i)}}.
\end{equation}

\noindent Subtracting the $i+1$ row from the $i$-th for all the rows of $G$ except for 
the last one nullifies all the elements of the first column except for the last element. 
Decomposing the determinant with respect to this column and factoring out 
common multiples, we see that 
$$
  det\(G(\tau)\) = \(1!\, 2!\, \cdots (m-1)!\, \tau^{m(m-1)/2}
  \prod\limits_{k=2}^m (b_k - b_1)\, \prod\limits_{2 \le i < j \le m} (b_j - b_i)^2\)^{-1} 
  det\({\tilde G}\),
$$

\noindent where ${\tilde G}$ is the following $(m-1) \times (m-1)$ matrix:
\begin{equation} \label{equation:Gtilde}
  {\tilde G}_{ij} = \prod\limits_{k=1}^{m-1-i} (b_{j+1} - b_1 -k\tau).
\end{equation}

\noindent For example, for $m = 4$,
$$
{\tilde G} = 
\[
  \begin{array}{ccc}
    (b_2 - b_1 - \tau)(b_2 - b_1 - 2\tau) & (b_3 - b_1 - \tau)(b_3 - b_1 - 2\tau) & 
    (b_4 - b_1 - \tau)(b_4 - b_1 - 2\tau) \\
     & & \\
    b_2 - b_1 - \tau & b_3 - b_1 - \tau & b_4 - b_1 - \tau \\
     & & \\
    1 & 1 & 1 \\
  \end{array}
\].
$$

\noindent To prove (\ref{equation:det(G)}), we have to prove that
$$
det\({\tilde G}\) = \prod\limits_{2 \le i < j \le m} (b_j - b_i).
$$

	Let us show that $det\({\tilde G}\)$ can be turned into the standard Vandermonde
without changing the determinant by means of row operations. First, let us take a look 
at our example. Multiplying the last row by $b_1 + \tau$ and adding the result to 
the previous row gives us the following matrix:
$$
\[
  \begin{array}{ccc}
    (b_2 - b_1 - \tau)(b_2 - b_1 - 2\tau) & (b_3 - b_1 - \tau)(b_3 - b_1 - 2\tau) & 
    (b_4 - b_1 - \tau)(b_4 - b_1 - 2\tau) \\
     & & \\
    b_2 & b_3 & b_4 \\
     & & \\
    1 & 1 & 1 \\
  \end{array}
\].
$$

\noindent Now let us take the last row times $-(b_1 + \tau)(b_1 + 2\tau)$ plus 
the second row times $2b_1 + 3\tau$ and add them to the first row. The result is
 the desired Vandermonde matrix:
$$
\[
  \begin{array}{ccc}
    b_2^2 & b_3^2 & b_4^2 \\
     & & \\
    b_2 & b_3 & b_4 \\
     & & \\
    1 & 1 & 1 \\
  \end{array}
\].
$$  

\noindent It is easy to complete the argument using induction. \\

	Let $J$ be the matrix with ones on the main diagonal and right above it, and 
zeros elsewhere:
$$
J_{ij} = \left\{ 
                     \begin{array}{ll}
                        1, & if~~ i=j-1; \\
                        1, & if~~ i=j; \\
                        0, & otherwise.
                     \end{array}
            \right.
$$ 
\noindent We want to prove that 
\begin{equation}
\label{equation:GconjX}
  G(\tau)\, EX(\tau)\, G(\tau)^{-1} = J.
\end{equation}

\noindent Consider the following $m \times m$ matrix:
$$
J'_{ij} = \left\{ 
                     \begin{array}{ll}
                        -1, & if~~ i=j-1; \\
                        1, & if~~ i=j; \\
                        0, & otherwise.
                     \end{array}
            \right.
$$ 

\noindent $J'$ is invertible, so $G(\tau)\, X(\tau) = J\, G(\tau)$ 
if and only if $J'\, G(\tau)\, X(\tau) = J'\, J\, G(\tau)$. It is not hard to see 
that 
$$
\(J' G(\tau)\)_{ij} = \left\{ 
\begin{array}{cl}
  0, & if~~ j = 1~~ and~~ i \ne m, \\
   & \\
  \displaystyle{\frac{1}{\prod\limits_{k=2}^m (b_1 - b_k)}}, & if~~ j = 1~~ and~~ i = m, \\
   & \\
  \displaystyle{\frac{\prod\limits_{k=0}^{m-i-1} (b_j - b_1 -k\tau)}
  {(m-i)!\, \tau^{m-i}\, \prod\limits_{k=1 \atop{k \ne j}}^m (b_j - b_k)}}, & if~~ j > 1
\end{array} \right.
$$

\noindent and that
$$
\(J' J\, G(\tau)\)_{ij} = \left\{ 
\begin{array}{cl}
  0, & if~~ j = 1~~ and~~ i < m-1, \\
   & \\
  \displaystyle{\frac{1}{\prod\limits_{k=2}^m (b_1 - b_k)}}, & if~~ j = 1~~ and~~ i = m - 1, \\
   & \\
  \displaystyle{\frac{1}{\prod\limits_{k=1 \atop{k \ne j}}^m (b_j - b_k)}}, & if~~ i = m, \\
   & \\
  \displaystyle{\frac{(b_j - b_1 + \tau)\, \prod\limits_{k=0}^{m-i-2} (b_j - b_1 -k\tau)}
  {(m-i)!\, \tau^{m-i}\, \prod\limits_{k=1 \atop{k \ne j}}^m (b_j - b_k)}}, & 
  if~~ j > 1~~ and~~ i < m.
\end{array} \right.
$$

	Let us prove that 
\begin{equation}
  \label{equation:J-J'}
  \sum\limits_{l=1}^m \(J'G(\tau)\)_{il} EX(\tau)_{lj} = 
  \(J' J\, G(\tau)\)_{ij}.  
\end{equation}

\noindent In the first case $j = 1$, $i = 1, \cdots ,m-2$, after factoring out and cancellation, 
(\ref{equation:J-J'}) becomes the following identity:
$$
\sum\limits_{l=2}^m \displaystyle{\frac{\prod\limits_{k=2}^{m-i-1} (b_l - b_1 - k\tau)}
{\prod\limits_{k=2 \atop{k \ne l}}^m (b_l - b_k)}} = 0.
$$

\noindent Up to a change of variables, this is identity (7.59) from \cite{G1}. \\

	In the second case $j = 1$, $i = m-1$, after factoring out and cancellation, 
(\ref{equation:J-J'}) becomes the following identity:
\begin{equation} \label{equation:idJ-J'1}
  \sum\limits_{l=2}^m \displaystyle{\frac{\prod\limits_{k=2 \atop{k \ne l}}^m 
  (b_1 - b_k + \tau)}{\prod\limits_{k=2 \atop{k \ne l}}^m (b_l - b_k)}} =1.
\end{equation}

\noindent Let us employ the strategy introduced at the beginning of the proof 
to  show that the left hand side of (\ref{equation:idJ-J'1}) is a polynomial in 
$b_i$ and $\tau$ of degree zero and thus a constant. Then set $b_1 + \tau = b_2$ 
to see that the constant is in fact $1$. \\

	In the third case $i = m$, $j = 1, \cdots ,m$, after some simplification, 
(\ref{equation:J-J'}) becomes the identity
$$
  \sum\limits_{l=1}^m \displaystyle{\frac{\prod\limits_{k=1 \atop{k \ne l}}^m 
  (b_j - b_k + \tau)}{\prod\limits_{k=1 \atop{k \ne l}}^m (b_l - b_k)}} =1,
$$

\noindent which is proven similarly to (\ref{equation:idJ-J'1}). \\

	Finally, for $i = 1, \cdots ,m-1$ and $j = 2, \cdots ,m$, (\ref{equation:J-J'}) 
boils down to the following identity:
$$
  \sum\limits_{l=2}^m \displaystyle{\frac
  {\prod\limits_{k=1}^{m-i-1} (b_l - b_1 - k\tau)}
  {\prod\limits_{k=0}^{m-i-2} (b_j - b_1 - k\tau)}}\, 
  \displaystyle{\frac{\prod\limits_{k=2 \atop{k \ne l}}^m (b_j - b_k + \tau)}
  {\prod\limits_{k=2 \atop{k \ne l}}^m (b_l - b_k)}} = 1.
$$

\noindent Our strategy works here again aided at some point by the identity (7.59) 
from \cite{G1}. \\

\noindent $\bullet$~~~ The last statement of the lemma follows from the last two
statements of Lemma \ref{lemma:Z}. Another way to see it is to observe 
that the first column of $G^{-1}(\tau)$ is $(m-1)!\, \tau^{m-1} e$. $\Box$ \\

	Let $V$ and $W$ be the matrices  composed of the eigenvectors $v_i$ 
(\ref{equation:v}) of the residue matrix $B$ and of the eigenvectors 
$w_i$ (\ref{equation:w}) of the residue matrix $C$ as columns respectively. 
For $\tau \in \C$, let $Z(\tau)$ be the following $m \times m$ matrix: 
\begin{equation} \label{equation:Z}
  Z_{ij}(\tau) = \displaystyle{\frac{\nu_j^2(\tau)}{b_j - c_i + a_2 + \tau}}.
\end{equation}

\begin{lemma} $\phantom{trr}$ \label{lemma:Z} \\ \\
  \begin{tabular}{ll}
  1. & $Z(0) = W^{-1}V$ \\
      & \\
  2. &  $Z^{-1}_{ij}(\tau) = \displaystyle{\frac{\mu_j^2(\tau)}{b_i - c_j + a_2 + \tau}}$ \\
      & \\
  3. & $Z^{-1}(\tau_1)\, Z(\tau_2) = EX(\tau_2 - \tau_1)$,~~~
         $Z(\tau_1)\, Z^{-1}(\tau_2) = EY(\tau_2 - \tau_1)$ \\
      & \\
  4. & $Z(\tau)\, e = e$
  \end{tabular}
\end{lemma} 

\proof Let us first prove that
\begin{equation} \label{equation:W0inv}
  W^{-1}_{ij} = \left\{
  \begin{array}{ll} 
  \displaystyle{\frac{1}{b_j - c_i + a_2}}\, 
  \prod\limits_{k=1}^{j-1} \displaystyle{\frac{c_k - c_i}{b_k - c_i + a_2}}, & if~~~ i \ge j, \\
   & \\
  0, & if~~~ i < j. \end{array} \right.
\end{equation}

\noindent To prove (\ref{equation:W0inv}), we have to show that
\begin{equation} \label{equation:W0inv:aux}
  \sum\limits_{l=1}^m W_{il}\, W_{lj}^{-1} = \delta_{ij}.
\end{equation}

\noindent Since both $W$ and $W^{-1}$ are lower triangular, it is clear that 
the right hand side of (\ref{equation:W0inv:aux}) equals $0$ for $j>i$. 
For $i=j$, the the right hand side of (\ref{equation:W0inv:aux}) has only one 
summand. It is easy to see that the summand equals $1$. Finally, for $i>j$ 
(\ref{equation:W0inv:aux}) boils down to the identity
$$
\sum\limits_{l=j}^i \displaystyle{\frac{\prod\limits_{k=j+1}^{i-1} (b_k - c_l + a_2)}
{\prod\limits_{k=j \atop{k \ne l}}^i (c_k - c_l)}} = 0,
$$
which is up to a change of variables equivalent to (7.59) from \cite{G1}. \\

\noindent 1.~~~ We need to show that 
$$
\sum\limits_{l=1}^m W^{-1}_{il} V_{lj} = Z(0)_{ij}
$$
or, equivalently,
\begin{equation} \label{equation:WinvV=Z(0)}
  \sum\limits_{l=1}^{\min(i,j)} \displaystyle{\frac{\prod\limits_{k=1 \atop{k \ne j}}^{l-1} 
  (b_j - b_k)}{\prod\limits_{k=1}^l (b_j - c_k + a_2)}} (b_j - c_i + a_2) 
  \displaystyle{\frac{\prod\limits_{k=1}^{l-1} (c_k - c_i)}{\prod\limits_{k=1}^l (b_k - c_i + a_2)}} 
  (b_l - c_l + a_2) = 1.
\end{equation}

\noindent Let us add up the summands of (\ref{equation:WinvV=Z(0)}) 
starting from the end. It is not hard to prove by induction that the sum of (the last) 
$n$ terms equals
$$
S_n = \prod\limits_{k=1}^{\min(i,j)-n}  
\displaystyle{\frac{(b_j - b_k)(c_k - c_i)}{(b_j - c_k + a_2)(b_k - c_i + a_2)}}.
$$ 

\noindent In particular, $S_{\min(i,j)} = 1$. \\

\noindent 2.~~~ The second statement of the lemma is equivalent to the rational 
Cauchy identity (\ref{equation:D_r^-1}). \\

\noindent 3.~~~ Let us prove that $\sum_{l=1}^m Z^{-1}(\tau_1)_{il}\, 
Z(\tau_2)_{lj} = EX(\tau_2 - \tau_1)_{ij}$. This is equivalent to proving 
the following identity:
$$
  \sum\limits_{l=1}^m \prod\limits_{k=1 \atop{k \ne i}}^m 
  \displaystyle{\frac{b_k - c_l + a_2 + \tau_1}{b_j - b_k + \tau_2 - \tau_1}}\, 
  \prod\limits_{k=1 \atop{k \ne l}}^m 
  \displaystyle{\frac{b_j - c_k + a_2 + \tau_2}{c_k - c_l}} = 1.
$$ 

\noindent The strategy of Lemma \ref{lemma:EX,EY} works here again aided by 
identity (7.59) from \cite{G1}. The corresponding identity for $EY$ is proven similarly. \\

\noindent 4.~~~ We need to show that for any $i = 1, \cdots ,m$,
$$
  \sum\limits_{j=1}^m \displaystyle
  {\frac{\prod\limits_{k=1 \atop{k \ne i}}^m (b_j - c_k + \tau)}
  {\prod\limits_{k=1 \atop{k \ne j}}^m (b_j - b_k)}} = 1.
$$

\noindent But this, up to a change of variables, is the identity (7.60) from \cite{G1}. 
$\Box$

\begin{lemma} \label{lemma:expX,Y}
$$
EX(\tau) = e^{X \tau} \mbox{~~and~~} EY(\tau) = e^{Y \tau}.
$$
\end{lemma}

\proof We shall only give a proof for $EX(\tau)$. According to the first statement of 
Lemma \ref{lemma:EX,EY}, $EX(\tau)$ is a one-parameter group. It is easy to see for an 
off-diagonal term that 
$$
\displaystyle{\frac{d}{d \tau}} \left. EX(\tau)_{ij} \right|_{\tau = 0} = 
\displaystyle{\frac{1}{b_j - b_i}}. 
$$
According to the second statement of Lemma \ref{lemma:EX,EY}, the Jordan normal 
form of the derivative at zero is a single block of full size $m$ with the eigenvalue $0$. 
According to the third statement of Lemma \ref{lemma:EX,EY}, the corresponding 
eigenvector is $e$. Thus, the sum of the elements of every row must add up to zero. $\Box$ \\

{\it Proof of Lemma \ref{lemma:propXY} --- } The first statement of the lemma is already 
proven in the course of proving Lemma \ref{lemma:expX,Y}. The second is proven by an 
easy computation. To prove the third statement of the lemma, let us prove that 
\begin{equation} \label{equation:VEXVinv<->WE-YWinv}
 V\, e^{X \tau}\, V^{-1} = W\, e^{-Y \tau}\, W^{-1}.
\end{equation} 
Rewriting this formula as 
$$
Z(0)\, EX(\tau) = EY(-\tau)\, Z(0)
$$
and using the third property of Lemma \ref{lemma:Z} proves 
(\ref{equation:VEXVinv<->WE-YWinv}). Differentiating (\ref{equation:VEXVinv<->WE-YWinv}) 
at zero finishes the proof for the bases $v_i$ (\ref{equation:v}) and $w_i$ (\ref{equation:w}). 
Recalling that the additive hypergeometric triple is rigid completes the argument. $\Box$ \\

{\it Proof of Lemma \ref{lemma:ev} --- } Let $B_d(0)$ and $C_d(0)$ be the Jordan normal 
forms $diag(b_1, \dots ,b_m)$ and $diag(-c_1, \dots ,-c_m)$ of the residue operators 
$B$ and $C$ respectively. Let $B_d(\tau_1) = B_d(0) + k_1 \tau_1 Id$ and 
$C_d(\tau_2) = C_d(0) + k_2 \tau_2 Id$. Recalling that $S = V X V^{-1} = -W Y W^{-1}$ 
and using the second statement of Lemma \ref{lemma:propXY}, it is not hard to show that 
the formulae
\begin{equation} \label{equation:flows}
  B(\tau_1) = V e^{X \tau_1} B_d(\tau_1) e^{-X \tau_1} V^{-1} \mbox{~~and~~} 
  C(\tau_2) = W e^{Y \tau_2} C_d(\tau_2) e^{-Y \tau_2} W^{-1}
\end{equation}

\noindent solve the dynamical equations (\ref{equation:main}), which is equivalent 
to the statements of Lemma \ref{lemma:ev}. $\Box$ \\ 

{\it Proof of Theorem \ref{thm:Frob_series} --- } Let us first prove the formula for $\(T_0\)_i$. 
The fact $S = V X V^{-1}$ implies that 
$$
  e^{nS} v_i = V e^{nX} e_i = \sum\limits_{j=0}^m EX(n)_{ji}\, v_j = 
  \sum\limits_{j=0}^m \xi_i^2(n)\, \displaystyle{\frac{v_j}{b_i - b_j + n}}. 
$$

\noindent Thus, we want to prove that 
\begin{equation} \label{equation:series:1}
  (T_0)_i = z^{b_i} \(v_i + \sum_{n=1}^{\infty} \alpha_{in}\, \xi_i^2(n)\, z^n\, 
  \sum\limits_{j=0}^m \displaystyle{\frac{v_j}{b_i - b_j + n}}\). 
\end{equation}

\noindent Let us rewrite the mHGS (\ref{equation:m-hgs}) as 
\begin{equation} \label{equation:mHGS:no_denominator}
  z \frac{d\,f}{d\,z}=\[ B - A(z+z^2+z^3+\cdots) \] f(z).
\end{equation}

\noindent Let us seek solutions to (\ref{equation:mHGS:no_denominator}) in the form 
\begin{equation} \label{equation:series:2}
  (T_0)_i = z^{b_i} \(\sum_{n=0}^{\infty} (T_0)_{in}\, z^n\). 
\end{equation}

\noindent Plugging (\ref{equation:series:2}) into (\ref{equation:mHGS:no_denominator}), 
we obtain the following recursive relation on $(T_0)_{in}$:
\begin{equation} \label{equation:recursion_at_0}
  \(T_0\)_{in} = \(B-(b_i + n) Id\)^{-1} A\, \((T_0)_{i0} + (T_0)_{i1} + \dots +(T_0)_{i,n-1}\),
\end{equation}

\noindent where $(T_0)_{i0}=v_i$. Let us prove that 
\begin{equation} \label{equation:recursion_at_0:1}
  \(T_0\)_{in} = \alpha_{in}\, \xi_i^2(n)\, \sum\limits_{j=0}^m \displaystyle{\frac{v_j}{b_i - b_j + n}}
\end{equation}

\noindent satisfies the recursion (\ref{equation:recursion_at_0}) by induction. First, let us 
establish the base. According to (\ref{equation:Ax}), $A v_i = a_2 v_i - (v_i, u)_r\, u$. For 
the case in consideration, $a_2 = 0$. Recalling that $u = v_1 + \dots + v_m$ and that 
the $v_i$ form an orthogonal basis with respect to the product (\ref{equation:(v,v),(w,w)}), 
we get $$A v_i = -\nu_{i\, r}^2\, \sum\limits_{j=1}^m v_j.$$ Applying $\(B-(b_i + 1) Id\)^{-1}$ 
to the right hand side of this formula, we obtain that 
$$
  (T_0)_{i1} = \nu_{i\, r}^2\, \sum\limits_{j=1}^m \displaystyle{\frac{v_j}{b_i - b_j +1}}.
$$
Observing that $\alpha_{i1}\, \xi_i^2(1) = \nu_{i\, r}^2$ finishes establishing the base. \\

	To make the step of induction from $n$ to $n+1$, we go along the same lines. In the end, 
the proof boils down to the identity 
$$
  \prod\limits_{k=1}^m \displaystyle{\frac{b_j - c_k + \tau}{b_j - b_k + \tau}} - 
  \sum\limits_{i=1}^m \displaystyle{\frac{1}{b_j - b_i + \tau}}\,
  \displaystyle{\frac{\prod\limits_{k=1}^m (b_i - c_k)}
  {\prod\limits_{k=1 \atop{k \ne i}}^m (b_i - b_k)}} = 1 
$$

\noindent (for $\tau = n$), which is proven using the method of Lemma \ref{lemma:expX,Y}. 
The proof for $\(T_{\infty}\)_i$ is similar. Finally, convergence follows from the following 
theorem: $\Box$ \\

\begin{theorem} \label{thm:hg:T_0,Tinf} 
  \begin{equation} \label{equation:back:T_0} (T_0)^j_i(z)=\end{equation}
    $$
      \left\{
      \begin{array}{l}
      z^{b_i} v_i^j\, \phantom{}_mF_{m-1}\(\left.
      \displaystyle{{b_i - c_1, \dots ,b_i - c_j, b_i - c_{j+1} + 1, \dots ,
      b_i - c_m + 1} \atop{b_i - b_1, \dots ,b_i - b_{j-1}, b_i - b_j + 1,\dots
      ,\widehat{b_i - b_i + 1}, \dots ,b_i - b_m + 1}} \right|\, z\), \mbox{~if~} i \ge j, \\
        \\
      \begin{array}{l}
      z^{b_i+1}\, (b_j - c_j)\, \nu_i^2(0)\, \displaystyle{\frac{\prod\limits_{k=j+1}^m
      (b_i - c_k + 1)}{\prod\limits_{k=j}^m (b_i - b_k + 1)}} \times 
      \phantom{there lived a fellow fro wealing who ha} \mbox{if~~} i < j; \\ 
        \\
      \phantom{}_mF_{m-1}\(\left. \displaystyle{{b_i - c_1 + 1, \dots
      ,b_i - c_j + 1, b_i - c_{j+1} + 2, \dots ,b_i - c_m + 2} \atop{b_i - b_1 + 1,\dots
      \widehat{b_i - b_i + 1} \dots ,b_i - b_{j-1} + 1, b_i - b_j + 2, \dots ,b_i - b_m + 2}} \right|\, z \),
      \end{array}
      \end{array} \right.
    $$ 
    \medskip 

  \begin{equation} \label{equation:back:Tinf} (T_{\infty})^j_i 
  \(\displaystyle{\frac1z}\) = \end{equation}
  $$
      \left\{
      \begin{array}{l}
      \begin{array}{l}
      z^{c_i-1}\, (b_j - c_j)\, \mu_i^2(0)\, \displaystyle{\frac{\prod\limits_{k=1}^{j-1} (b_k - c_i + 1)}
      {\prod\limits_{k=1}^j (c_k - c_i + 1)}} \times 
      \phantom{there lived a young fellow} \mbox{if~} i > j,\\ \\
      \phantom{}_mF_{m-1}\(\left. \displaystyle{{b_1 - c_i + 2,\dots
      ,b_{j-1} - c_i + 2, b_j - c_i + 1, \dots ,b_m - c_i + 1}\atop{c_1 - c_i + 2, \dots ,c_j - c_i + 2,
      c_{j+1} - c_i + 1, \dots ,\widehat{c_i - c_i + 1} \dots ,c_m - c_i + 1}} \right|\, 
      \displaystyle{\frac1z} \),
      \end{array} \\
       \\ \\
      z^{c_i}\, w_i^j\, \phantom{}_mF_{m-1}\(\left.
      \displaystyle{{b_1 - c_i + 1, \dots ,b_{j-1} - c_i + 1, b_j - c_i, \dots ,b_m - c_i}
      \atop{c_1 - c_i + 1, \dots , \widehat{c_i - c_i + 1}, \dots ,c_j - c_i + 1, c_{j+1} - c_i,
      \dots ,c_m - c_i}} \right|\, \displaystyle{\frac1z} \), \\ \\
      \mbox{if~} i \le j.
      \end{array} \right.
    $$
\end{theorem}

	{\it Proof--- ~~} According to (\ref{equation:recursion_at_0:1}), 
$$
  \(T_0\)_{in} = \alpha_{in}\, \xi_i^2(n)\, \sum\limits_{j=0}^m \displaystyle{\frac{v_j}{b_i - b_j + n}}. 
$$
Let us call $\(T_0\)_n$ (not to be confused with $(T_0)_i$ from Theorem \ref{thm:Frob_series}) 
the matrix composed of columns $\(T_0\)_{in}$ and let us call ${\cal A}_n$ 
the diagonal $m \times m$ matrix with $\alpha_{in}$ on the diagonal. Then
\begin{equation} \label{equation:T_0:v}
 \(T_0\)_n = V e^{nX} {\cal A}_n = W W^{-1} V e^{nX} {\cal A}_n.
\end{equation}
According to the first statement of Lemma \ref{lemma:Z}, $W^{-1} V = Z(0)$.
According to the third statement of Lemma \ref{lemma:Z} combined with 
Lemma \ref{lemma:expX,Y}, $e^{nX} = Z^{-1}(0) Z(n)$. Thus, we can rewrite 
(\ref{equation:T_0:v}) as
\begin{equation} \label{equation:T_0:w}
  \(T_0\)_n = W Z(n)\, {\cal A}_n = W D_r^t(n)\, {\cal N}^2(n)\,{\cal A}_n.
\end{equation}

\noindent Let us prove that
\begin{equation} \label{equation:WD}
  \(W D_r^t(n)\)_{ij} = (b_i - c_i)\, \displaystyle{\frac{\prod\limits_{k=1}^{i-1} (b_j - b_k + n)}
  {\prod\limits_{k=1}^i (b_j - c_k + n)}}.
\end{equation}

\noindent This boils down to the identity
$$
\sum\limits_{l=1}^i \displaystyle{\frac{\prod\limits_{k=1 \atop{k \ne l}}^i (b_j - c_k + n)}
{\prod\limits_{k=1}^{i-1} (b_j - b_k + n)}}\, 
\displaystyle{\frac{\prod\limits_{k=1}^{i-1} (b_k - c_l)}
{\prod\limits_{k=1 \atop{k \ne l}}^i (c_k - c_l)}} = 1,
$$

\noindent proven using the strategy of Lemma \ref{lemma:expX,Y}. 
A look at the definition (\ref{equation:HG:F:B-H}) of the generalized hypergeometric 
function finishes the proof of (\ref{equation:back:T_0}). \\

	Similarly,
$$
 \(T_{\infty}\)_n = W e^{nY} {\cal B}_n = V Z^{-1}(n)\, {\cal B}_n = 
  V D_r(n)\, {\cal M}^2(n)\, {\cal B}_n.
$$

\noindent Similarly to (\ref{equation:WD}),
$$
\(V D_r(n)\)_{ij} = (b_i - c_i)\, \displaystyle{\frac{\prod\limits_{k=i+1}^m (c_k - c_j + n)}
{\prod\limits_{k=i}^m (b_k - c_j + n)}}, 
$$
which proves (\ref{equation:back:Tinf}). $\Box$ \\

	The formula (\ref{equation:trig_continuation}) follows from Theorem \ref{thm:hg:T_0,Tinf} 
combined with (\ref{equation:analytic_cont_of_mFm-1_to_infty}). \\

	{\it Proof of Theorems \ref{thm:(F,F):ell} and \ref{thm:D_ell^-1} --- }
Recall that the space $\H$ has two special bases $\(F_0\)_i$ and $\(F_{\infty}\)^{\dagger}_j$, 
such that $\(\(F_0\)_i, \(F_{\infty}\)^{\dagger}_j\) = 1/sn(a_2 + b_i - c_j)$. Let us set
$\(\(F_0\)_i, \(F_0\)_j\) = \delta_{ij} / \nu_{i\, ell}^2$. Let $\(F_{\infty}\)^{\dagger}_j = 
\sum_{k=1}^m x_{kj} \(F_0\)_k$. Taking the scalar product of both sides of this equality 
with $\(F_0\)_i$, we see that 
$$
  \(F_{\infty}\)^{\dagger}_j = \sum\limits_{i=1}^m 
  \displaystyle{\frac{\nu_{i\, ell}^2}{sn(a_2 + b_i - c_j)}} \(F_0\)_i.
$$
\noindent Our construction has the following symmetry: if we switch the points 
$0$ and $\infty$ on the Riemann sphere and simultaneously switch the $b_i$ and $-c_i$ for 
all $i = 1, \dots ,m$, all the formulae remain valid. Applying the symmetry to the above formula, 
we get
$$
  \(F_0\)_j = \sum\limits_{i=1}^m 
  \displaystyle{\frac{\mu_{i\, ell}^2}{sn(a_2 + b_j - c_i)}} \(F_{\infty}\)^{\dagger}_i.
$$
\noindent Rewriting the last two formulae in the matrix form as
$$
F_{\infty}^{\dagger} = F_0 D_{ell}^t\, {\cal N}_{ell}^2 \mbox{~~and~~} 
F_0 = F_{\infty}^{\dagger} D_{ell}\, {\cal M}_{ell}^2
$$
and comparing them to each other proves Theorem \ref{thm:D_ell^-1}. \\

	The formula $ \({\cal N}\, D\, {\cal M}\)^t = \({\cal N}\, D\, {\cal M}\)^{-1}$ implies 
that $\(\(F_{\infty}\)^{\dagger}_i / \mu_{i\, ell}, \(F_{\infty}\)^{\dagger}_j / \mu_{j\, ell}\) = \delta_{ij}$. 
Finally, the dimension count proves the uniqueness of the product. $\Box$ \\

	{\it Proof of Theorem \ref{thm:hk} --- } Recall that in Theorem \ref{thm:hk} 
$\omega_1 = 1$, $\omega_2$ is an imaginary number, all the local exponents are real, 
and $0 \le b_1 < \dots < b_m < 1$ and $0 \le c_1 < \dots < c_m < 1$. In this case, 
all the values of $sn(a_2 + b_i - c_j)$ are real for $i,j = 1, \dots ,m$. If in addition the 
positivity conditions (\ref{equation:hg:form:pos}) are satisfied, then the real 
symmetric form (\ref{equation:(F,F):ell}) on $\H$ is sign-definite as well as its counterpart 
on $\H'$. They can be extended to complex numbers either in the hermitian or in the complex 
symmetric fashion. The quaternionic action (\ref{equation:quaternions}) applied to the fields 
(\ref{equation:generating_functions}) produces the action (\ref{equation:quaternions:F}). It is 
not hard to see that in the sign-definite case the action (\ref{equation:quaternions:F}) is 
hyperk\"ahler. $\Box$ \\

\noindent {\bf Theorem \ref{thm:product} -- extended version ~~~} \label{thm:product:ext} \\ \\ {\it
\noindent $\bullet~~$ For any  complex times $\tau_1$ and $\tau_2$, 
\begin{equation} \label{equation:v<->w}
  v_i(\tau_1) = \sum\limits_{j=1}^m \displaystyle{\frac{\nu_i^2(\tau_1 + \tau_2)}
  {a_2 + b_i - c_j + \tau_1 + \tau_2}}\, w_j(\tau_2).
\end{equation}

\noindent There exists a unique up to a constant multiple complex symmetric scalar 
product $(*,*)_{\tau_1, \tau_2}$ on $H_{\tau_1, \tau_2}$ such that the bases 
$v_i(\tau_1)$ and $w_i(\tau_2)$ are simultaneously orthogonal with respect to it:
  $$ 
     \(v_i(\tau_1), v_j(\tau_1)\)_{\tau_1, \tau_2} = 
     \delta_{ij}\, \nu_i^2\(\tau_1 + \tau_2\), ~~~
     \(w_i(\tau_2), w_j(\tau_2)\)_{\tau_1, \tau_2} = 
     \delta_{ij}\, \mu_i^2\(\tau_1 + \tau_2\).  
  $$ 
  
\noindent (Note that formula (\ref{equation:(v,w)}) of the original Theorem \ref{thm:product} 
follows from (\ref{equation:v<->w}) and (\ref{equation:(*,*)_tau12}) combined.) \\

\noindent $\bullet~~$ For the bases $v_i(\tau_1)$ and $v_i(\tau_2)$ such that $\tau_1 \ne \tau_2$, 
\begin{equation} \label{equation:transition:+}
  v_j(\tau_1) = \sum\limits_{i=1}^m \displaystyle{\frac{\xi_j^2(\tau_1 -  \tau_2)}
  {b_j - b_i + \tau_1 - \tau_2}}\, v_i(\tau_2).
\end{equation}

\noindent There exists a unique up to a constant multiple complex symmetric scalar 
product $(*,*)_{\tau_1, \tau_2}^+$ such that the bases $v_i(\tau_1)$ and $v_i(\tau_2)$ 
are simultaneously orthogonal with respect to it: 
\begin{equation} \label{equation:(*,*)^+}  
  \(v_i(\tau_1), v_j(\tau_1)\)_{\tau_1, \tau_2}^+ = \delta_{ij}\, \xi_i^2(\tau_1 - \tau_2), ~~~ 
  \(v_i(\tau_2), v_j(\tau_2)\)_{\tau_1, \tau_2}^+ = -\delta_{ij}\, \xi_i^2(\tau_2 - \tau_1).
\end{equation}

\noindent $\bullet~~$ For the bases $w_i(\tau_1)$ and $w_i(\tau_2)$ such that $\tau_1 \ne \tau_2$, 
\begin{equation} \label{equation:transition:-}
  w_j(\tau_1) = \sum\limits_{i=1}^m \displaystyle{\frac{\theta_j^2(\tau_1 -  \tau_2)}
  {c_i - c_j + \tau_1 - \tau_2}}\, w_i(\tau_2).
\end{equation}

\noindent There exists a unique up to a constant multiple complex symmetric scalar 
product $(*,*)_{\tau_1, \tau_2}^-$ such that the bases $w_i(\tau_1)$ and $w_i(\tau_2)$ 
are simultaneously orthogonal with respect to it: 
\begin{equation} \label{equation:(*,*)^-}  
  \(w_i(\tau_1), w_j(\tau_1)\)_{\tau_1, \tau_2}^- = \delta_{ij}\, \theta_i^2(\tau_1 - \tau_2), ~~~ 
  \(w_i(\tau_2), w_j(\tau_2)\)_{\tau_1, \tau_2}^- = -\delta_{ij}\, \theta_i^2(\tau_2 - \tau_1).
\end{equation} }

\proof $~~~\bullet~~$ The formula (\ref{equation:v<->w}) in the matrix form reads as 
$V e^{X \tau_1} = W e^{Y \tau_2} Z(\tau_1 + \tau_2)$.  
By Lemma \ref{lemma:Z}, $e^{Y \tau_2} = Z(\tau_1) Z^{-1}(\tau_1 + \tau_2)$, so we get
$V e^{X \tau_1} = W Z(\tau_1)$. 
By the same lemma, $V^{-1} W = Z^{-1}(0)$ and $Z^{-1}(0) Z(\tau_1) = e^{X \tau_1}$. 
The existence of the product (\ref{equation:(*,*)_tau12}) follows from the rational Cauchy 
identity (\ref{equation:D_r^-1}). Its uniqueness follows from the dimension count. \\

\noindent $\bullet~~$ The formula (\ref{equation:transition:+}) in matrix notations is simply
$V e^{X \tau_1} = V e^{X \tau_2} e^{X (\tau_1 - \tau_2)}$. For $\tau \ne 0$, let us set in the 
rational Cauchy identity $a_2 = 0$ and $c_i = b_i$ for $i = 1, \dots ,m$. The special case 
$$
  \[\displaystyle{\frac{1}{b_i - b_j + \tau}}\]^{-1}_{ij} = 
  -\displaystyle{\frac{\xi_i^2(-\tau)\, \xi_j^2(\tau)}{b_j - b_i + \tau}}
$$
\noindent gives rise to (\ref{equation:(*,*)^+}). Uniqueness follows from the dimension count. 
The last case is proven similarly. $\Box$ \\

	{\it Proof of Lemma \ref{lemma:A} --- } Let us begin with proving the third statement of 
the lemma first. To prove that $B(\tau_1) x + C(\tau_2) x + a_2(\tau_1, \tau_2) x = 
(x, u)_{\tau_1, \tau_2} u$ for any $x \in H_{\tau_1, \tau_2}$, it is enough to prove that this 
formula holds for any basis vector $v_i(\tau_1)$. The formula (\ref{equation:v<->w}) 
gives us
$$
  \(B(\tau_1) + a_2(\tau_1, \tau_2)\) v_i(\tau_1) = (b_i(\tau_1) + a_2(\tau_1, \tau_2)) 
  \sum\limits_{j=1}^m \displaystyle{\frac{\nu_i^2(\tau_1 + \tau_2)}
  {a_2 + b_i - c_j + \tau_1 + \tau_2}}\, w_j(\tau_2).
$$
\noindent On the other hand, 
$$
  C(\tau_2) v_i(\tau_1) = \sum\limits_{j=1}^m \displaystyle{\frac{\nu_i^2(\tau_1 + \tau_2)\, 
  c_j(\tau_2)}{a_2 + b_i - c_j + \tau_1 + \tau_2}}\, w_j(\tau_2). 
$$
\noindent But $b_i(\tau_1) + c_j(\tau_2) + a_2(\tau_1, \tau_2) = a_2 + b_i - c_j + \tau_1 + \tau_2$. 
Thus, 
$$
  \(B(\tau_1) + C(\tau_2) + a_2(\tau_1, \tau_2)\) v_i(\tau_1) = \nu_i^2(\tau_1 + \tau_2)\, 
  \sum\limits_{j=1}^m w_j(\tau_2) = (v_i(\tau_1), u)_{\tau_1, \tau_2} u. 
$$

	The first statement easily follows from the third. To prove the forth statement, we need 
to show that $\sum_{i=1}^m \nu_i^2(\tau_1 + \tau_2) = m a_2(\tau_1, \tau_2) + 
\sum_{i=1}^m (b(\tau_1) + c(\tau_2)) = \sum_{i=1}^m (a_2 + b_i - c_j + \tau_1 + \tau_2)$.
This is the formula (7.61) from \cite{G1} in different notations. Finally, the second statement 
follows from the third and the forth combined. $\Box$ 

\subsection{Beukers and Heckman revisited} \label{subsec:BH}%%%%%%%%%%%%%%
The main purpose of this subsection is to prove Lemma \ref{lemma:p<->q}. 
As a side product, we give a proof to Theorem \ref{thm:(p,p),(q,q)} different from those 
of Beukers and Heckman and of Haraoka. \\

	Recall that $p_i$ are the eigenvectors of the monodromy operator $M_0$ corresponding 
to the eigenvalues $e^{2 \pi \sqrt{-1}\, b_i}$, $q_i$ are the eigenvectors of the monodromy 
operator $M_{\infty}$ corresponding to the eigenvalues $e^{-2 \pi \sqrt{-1}\, c_i}$, and $r$ is 
the eigenvector of the monodromy operator $M_1$ corresponding to the eigenvalue 
$e^{2 \pi \sqrt{-1}\, a_1}$. Consider the following three flags in the space of solutions: 
${\cal F}_0 = Span(p_1) \subset Span(p_1, p_2) \subset \dots \subset \C^m$,  
${\cal F}_{\infty} = Span(q_m) \subset Span(q_m, q_{m-1}) \subset \dots \subset \C^m$, and 
${\cal F}_1 = Span(r) \subset \C^m$. The first two flags are complete. For two complete flags, 
one can always choose a basis $e_1, \cdots ,e_m$ making the flags opposite: 
${\cal F}_0 = Span(e_1) \subset Span(e_1, e_2) \subset \dots \subset \C^m$ and 
${\cal F}_{\infty} = Span(e_m) \subset Span(e_m, e_{m-1}) \subset \dots \subset \C^m$. 
In such a basis, the matrix of $M_0$ will look upper- and the matrix of $M_{\infty}$ 
lower-triangular.

\begin{theorem} \label{thm:M}
Let $b_1, \cdots , b_m$ and $c_1, \cdots ,c_m$ be generic complex numbers and let
\begin{equation} \label{equation:hg:M_0}
    (M_0)_{ij} = \left\{
    \begin{array}{lll}
        0, & if & i>j \\
        e^{2 \pi \sqrt{-1}\, b_i}, & if & i=j \\
        e^{2 \pi \sqrt{-1}\, \((j - i -1)a_2 + b_j + \sum_{k = i+1}^{j-1} (b_k - c_k)\)} 
        \(e^{2 \pi \sqrt{-1}\, (b_i - c_i + a_2)}-1\), & if & i<j \\
    \end{array} \right. , 
\end{equation}

\begin{equation} \label{equation:hg:M_1}
    (M_1 - e^{2 \pi \sqrt{-1}\, a_2} Id)_{ij} = 
    e^{-2 \pi \sqrt{-1}\, \((i - 1)a_2 + \sum_{k = 1}^i (b_k - c_k)\)} 
    \(1 - e^{2 \pi \sqrt{-1}\, (b_i - c_i + a_2)}\), 
\end{equation}

\begin{equation} \label{equation:hg:M_infty}
    (M_{\infty})_{ij} = \left\{
        \begin{array}{lll}
        e^{-2 \pi \sqrt{-1}\, (b_i + a_2)} \(e^{2 \pi \sqrt{-1}\, (b_i - c_i + a_2)}-1\), & if & i>j \\
        e^{-2 \pi \sqrt{-1}\, c_i}, & if & i=j \\
        0, & if & i<j \\
    \end{array} \right. .
\end{equation}

\noindent Then $M_{\infty} M_1 M_0 = Id$.
\end{theorem}

\proof Let us first prove that 
\begin{equation} \label{equation:M0^-1}
  \(M_0^{-1}\)_{ij} = \left\{ \begin{array}{lll}
    0, & if & i>j; \\
     & & \\
    e^{-2 \pi \sqrt{-1}\, b_i}, & if & i=j; \\
     & & \\
    e^{2 \pi \sqrt{-1}\, (a_2 - c_i)} \(e^{2 \pi \sqrt{-1}\, (c_i - b_i - a_2)} - 1\), & if & i<j.
  \end{array} \right.
\end{equation}

\noindent For $i > j$, $\sum_{k=1}^m \(M_0^{-1}\)_{ik} \(M_0\)_{kj} = 0$, since both $M_0$ 
and $M_0^{-1}$ are upper-triangular. For $i = j$, $\sum_{k=1}^m \(M_0^{-1}\)_{ik} \(M_0\)_{kj} = 
e^{-2 \pi \sqrt{-1}\, b_i}\, e^{2 \pi \sqrt{-1}\, b_i} = 1$. For $i < j$, 
$\sum_{k=1}^m \(M_0^{-1}\)_{ik} \(M_0\)_{kj} = \sum_{k=i}^j \(M_0^{-1}\)_{ik} \(M_0\)_{kj}$. 
After factoring out the common multiple $e^{2 \pi \sqrt{-1}\, (c_i - b_i - a_2)} - 1$, it is not hard 
to see that the sum telescopes to zero. The proof of the fact that $M_{\infty} M_1 = M_0^{-1}$ 
is also a straitforward computation with telescoping sums for three different cases $i > j$, $i = j$, 
and $i < j$. $\Box$ \\ 

	We shall suppress detailed proofs in the remaining part of this subsection. All the below 
formulae are not hard to prove by direct computation. After factoring out common multiples, 
simplification, and sometimes telescoping, all of them either boil down to the identities 
from \cite{G1} or become trivial just as above. \\

	For $i = 1, \dots ,m$, the vectors 

\begin{equation} \label{equation:p}
  p^j_i = \left\{ \begin{array}{clr}
    e^{2 \pi \sqrt{-1}\, (a_1 - a_2 + b_i - b_j)} \(e^{2 \pi \sqrt{-1}\, (b_j - c_j + a_2)} - 1\)\, \times & & \\
     & & if~~ i \ge j, \\ 
    \displaystyle{\frac{\prod\limits_{k = j+1}^m \(e^{2 \pi \sqrt{-1} (b_i - c_k + a_2)} - 1\)}
    {\prod\limits_{k = j \atop{k \ne i}}^m \(e^{2 \pi \sqrt{-1}\, (b_i - b_k)} - 1\)}}, & & \\ 
     & \phantom{trr} & \\
    0, & & if~~ i < j; \end{array} \right.
\end{equation}

\noindent are the eigenvectors of the matrix (\ref{equation:hg:M_0}) corresponding to
the eigenvalues $e^{2 \pi \sqrt{-1}\, b_i}$ and the vectors

\begin{equation} \label{equation:q}
  q^j_i = \left\{ \begin{array}{clr} 
    0, & & if~~ i > j; \\
     & & \\ & & \\
    e^{-2 \pi \sqrt{-1}\, \( j a_2 + \sum\limits_{k=1}^j (b_k - c_k)\)}\, 
    \(e^{2 \pi \sqrt{-1}\, (b_j - c_j + a_2)} - 1\) \times & \\ 
     & & if~~ i \le j. \\
    \displaystyle{\frac{\prod\limits_{k=1}^{j-1} \(e^{2 \pi \sqrt{-1}\, (b_k - c_i + a_2)} - 1\)}
    {\prod\limits_{k=1 \atop{k \ne i}}^j \(e^{2 \pi \sqrt{-1}\, (c_k - c_i)} - 1\)}}, & & 
  \end{array} \right.
\end{equation}

\noindent are the eigenvectors of the matrix (\ref{equation:hg:M_infty}) corresponding to 
the eigenvalues $e^{-2 \pi \sqrt{-1}\, c_i}$. The vector 

\begin{equation} \label{equation:r}
  r^i = e^{-2 \pi \sqrt{-1}\, \(ia_2 + \sum\limits_{k=1}^i (b_k - c_k)\)} 
  \(e^{2 \pi \sqrt{-1}\, (b_i - c_i + a_2)} - 1\)
\end{equation}

\noindent is the eigenvector of the matrix (\ref{equation:hg:M_1}) corresponding to the 
eigenvalue $e^{2 \pi \sqrt{-1}\, a_1}$ such that 
$$
  \sum\limits_{i = 1}^m p_i = \sum\limits_{i = 1}^m q_i = r.
$$

	Let $P$ and $Q$ be the $m \times m$ matrices composed of the vectors $p_i$ 
and $q_i$ as columns. Then

\begin{equation} \label{equation:PinvQ}
  \(P^{-1} Q\)_{ij} = \displaystyle{\frac{\prod\limits_{k = 1 \atop{k \ne i}}^m 
  \(e^{2 \pi \sqrt{-1}\, (b_k - c_j + a_2)} - 1\)}{\prod\limits_{k = 1 \atop{k \ne j}}^m 
  \(e^{2 \pi \sqrt{-1}\, (c_k - c_j)} - 1\)}} = 
  e^{-\pi \sqrt{-1}\, (a_1 + b_i - c_j)}\, \displaystyle{\frac{\prod\limits_{k = 1 \atop{k \ne i}}^m 
  \sin \pi (b_k - c_j + a_2)}{\prod\limits_{k = 1 \atop{k \ne j}}^m \sin \pi (c_k - c_j)}}
\end{equation}

\noindent or, in matrix notations, 

\begin{equation} \label{equation:PinvQ:matrix}
  P^{-1} Q = e^{-\pi \sqrt{-1}\, a_1}\, diag(e^{-\pi \sqrt{-1}\, b_i})\, D_{trig}\, {\cal M}_{trig}^2\, 
  diag(e^{\pi \sqrt{-1}\, c_j}),
\end{equation}
\smallskip

\noindent where $D_{trig}$ is the matrix (\ref{equation:D_trig}). Inverting 
(\ref{equation:PinvQ:matrix}) with the help of the trigonometric Cauchy identity 
(\ref{equation:D_trig^-1}) proves Lemma \ref{lemma:p<->q}. Theorem \ref{thm:(p,p),(q,q)} 
follows from Lemma \ref{lemma:p<->q} combined with the trigonometric Cauchy identity. 

%%%%%%%%%%%%%%%%%%%%%%%%%%%%%%%%%%%%%%%%%%%%%%%%
\section{Remarks and open questions}%%%%%%%%%%%%%%%%%%%%%%%%%%%%%

	A different and very interesting approach to the generalized hypergeometric function 
through fermionic fields was made recently in \cite{OrSch1} and \cite{OrSch2}. 
The following observation may serve as a starting point in investigating the relations 
between their view of the generalized hypergeometric function and ours. \\

Let $d(b)$ and $d(c)$ be $m \times m$ diagonal matrices with the $i$-th 
diagonal elements equal to
\begin{equation} \label{equation:d}
  \begin{array}{lll}
    \prod\limits_{k=1 \atop{k \ne i}}^m (b_i - b_k) & \mbox{and} &
    \prod\limits_{k=1 \atop{k \ne i}}^m (c_k - c_i)
  \end{array}
\end{equation}

\noindent respectively. Let ${\tilde X} = d^{-1}(b)\, X d(b)$ and 
${\tilde Y} = d^{-1}(c)\, Y d(c)$. Let us call $Vnd(x, \tau)$ the following
$m \times m$ Vandermonde matrix:
\begin{equation} \label{equation:Vnd}
  \[\begin{array}{cccc}
     1                           & 1                            & \cdots & 1                      \\
     x_1 + \tau             & x_2 + \tau              & \cdots & x_m + \tau        \\
     (x_1 + \tau)^2       & (x_2 + \tau)^2        & \cdots & (x_m + \tau)^2 \\
     \vdots                    & \vdots                    &           & \vdots               \\
     (x_1 + \tau)^{m-1} & (x_2 + \tau)^{m-1} & \cdots & (x_m + \tau)^{m-1}
  \end{array}\].
\end{equation}

\noindent Then $e^{\tau {\tilde X}} = Vnd^{-1}(b,0)\, Vnd(b,\tau)$ and
$e^{\tau {\tilde Y}} = Vnd^{-1}(c,0)\, Vnd(c,\tau)$. \\

	Also, it would be interesting to understand what unknown feature of our construction 
manifests itself through the existence of the ``extra'' products (\ref{equation:(*,*)^+}) and 
(\ref{equation:(*,*)^-}) in the extended version of Theorem \ref{thm:product}. 

\section{Appendix: rational Calogero-Moser system} \label{sec:CM}%%%%%%%%%%%%%%

	The following description of the rational Calogero-Moser system is taken from \cite{Wil}. 
Some original notations are changed in order to conform to those of this paper. See \cite{Wil}
for details and further references. \\

	For a positive integer $m$, let $Id$ be the $m \times m$ identity matrix. Let $V_m$ be 
the complex vector space of all quadruples $(B, X; v, w)$, where $B$ and $X$ are $m \times m$ matrices and $v$ and $w$ are column and raw vectors respectively of length $m$. 
Let $\tilde{C}_m$ be the subvariety of $V_m$ satisfying the equation
\begin{equation} \label{equation:rCM}
  [B,X] = v \otimes w - Id.
\end{equation}

\noindent Let the group $GL(m, \C)$ act on $V_m$ by
\begin{equation} \label{equation:GLactionVm}
  g\, \circ (B, X; v, w) = (gBg^{-1}, gXg^{-1}; gv, g^{-1} w).
\end{equation}

\noindent The action clearly preserves $\tilde{C}_m$. Let $C_m = \tilde{C}_m/GL(m, \C)$.
\begin{proposition} (\cite{Wil}) \label{prop:Wil}
  Let $(B, X; v, w) \in \tilde{C}_m$, and suppose $B$ is diagonalizable. Then the eigenvalues 
  of $B$ are distinct; and the $GL(m, \C)$-orbit of $(B, X; v, w)$ contains a representative 
  $B_d$ such that \\
  1.~~ $B_d$ is diagonal, $B_d = diag(b_1, \dots ,b_m)$; \\ 
  2.~~ all the entries of the vectors $v$ and $w$ are equal to $1$; \\
  3.~~ $X$ is a Calogero-Moser matrix, i.e. 
  $$
    X_{ij} = \left\{ \begin{array}{lll}
    p_i & \mbox{,if} & i = j, \\ & & \\
    (b_j - b_i)^{-1} & \mbox{,if} & i \ne j.
    \end{array}\right.
  $$
The representation is unique up to a simultaneous permutation of the parameters 
$p_i$ and $b_i$.
\end{proposition}

	Let us identify $gl(m, \C)$ and $\C^m$ with their dual spaces in the usual way:  
by means of the trace form on $gl(m, \C)$ and of the standard bilinear form on $\C^m$. 
Then we can consider $V_m$ as a cotangent bundle to $gl(m, \C) \oplus \C^m$, thinking 
of $(X,w)$ as a cotangent vector at the point $(B,v)$ of the space. Thus, $V_m$ is equipped 
with the complex symplectic form
\begin{equation} \label{w_on_Vm}
  \omega = \mbox{tr}\(dX \wedge dB + dw \wedge dv \)
\end{equation}

\noindent The map $\mu: V_m \to gl(m,\C)$
\begin{equation} \label{equation:mu} 
  \mu(B, X; v, w) = [B, X] - v \otimes w 
\end{equation}

\noindent is the moment map for the action (\ref{equation:GLactionVm}) such that 
$\tilde{C}_m = \mu^{-1}(-Id)$. The {\it Calogero-Moser flows} are the flows on $C_m$ 
induced by the $GL(m, \C)$-invariant flows
\begin{equation} \label{equation:C-M_flows}
  (B, X; v, w) \mapsto (B + k t_k X^{k-1}, X; v, w)
\end{equation}

\noindent on $\tilde{C}_m$. Let $C'_m$ be the subset of $C_m$ with $B$ diagonalizable. 
$C'_m$ is dense and open in $C_m$. The symplectic form induced on $C'_m$ takes the 
canonical form $\sum dp_i \wedge dx_i$, so $C'_m$ can be identified with the phase space 
of a system of $m$ classical indistinguishable particles in the complex plane with the 
Hamiltonian
\begin{equation} \label{equation:rCM}
 -H_2 = \frac12 \mbox{tr} X^2 = 
 \frac12 \sum\limits_{k = 1}^m p_k^2 - \sum\limits_{1 \le j < k \le m} (x_j - x_k)^{-1}, 
\end{equation} 

\noindent called the rational Calogero-Moser system. 

%%%%%%%%%%%%%%%%%%%%%%%%%%%%%%%%%%%%%%%%%%%%%%%%%
%%%%%%%%%%%%%%%%%%%%%%%%%%%%%%%%%%%%%%%%%%%%%%%%%
\section*{Acknowledgments}

    I would like to thank F.~Andrianov for many interesting conversations about Mathematics, 
its spirit and philosophy which have influenced this paper in a variety of ways. I would like to 
thank P.~Brosnan for explaining some algebraic geometry and especially for an important 
example of a rigid local system. I would like to thank X.~Wang for answering a lot of 
questions about the geometry of various objects mentioned in the paper. \\

	I would like to thank people who took their time to listen to the raw versions of 
the paper's results and/or to read the manuscript: A.~Borodin, P.~Etingof, M.~Hitrick, 
E.~D'Hoker, I.~Krichever, A.~Postnikov, E.~Rains, C.~Simpson, T.~Tao, V.~Varadarajan, 
and the referee of the paper. Their comments and questions were most helpful. \\ 

	The last, but not the least, I would like to thank the Master Class program of Utrecht 
University, the Netherlands, and my adviser then, G.~Heckman. They gave the author so much 
in one year, that it took him ten years to sort things out. In particular, this paper is a natural continuation of \cite{BH} lectured to us by F.~Beukers back in Utrecht. 

%%%%%%%%%%%%%%%%%%%%%%%%%%%%%%%%%%%%%%%%%%%%%%%%%%%
%%%%%%%%%%%%%%%%%%%%%%%%%%%%%%%%%%%%%%%%%%%%%%%%%%%

\medskip


\begin{thebibliography}{30}

\bibitem{BH} F.~Beukers, G.~Heckman, {\it Monodromy for the hypergeometric function
$\phantom{|}_nF_{n-1}$}, Inventiones Mathematicae 95 (1989), 325-354

% \bibitem{Bo} A.~Bolibrukh, {\it The 21st Hilbert problem for linear Fuchsian systems},
% proc. of the Steklov Institute of Mathematics, vol.206, AMS, Providence, RI, 1995

\bibitem{DR} M.~Dettweiler, S.~Reiter, {\it An Algorithm of Katz and its applications
to the inverse Galois problem}, J. Symbolic Computations (2000) 30, 761-798

\bibitem{G1} O.~Gleizer, {\it Explicit solutions of the additive Deligne-Simpson
problem and their applications}, Advances in Mathematics 178 (2003) 311-374

\bibitem{HC} A.~Hurwitz and R.~Courant, {\it Vorlesungen {\"u}ber allgemeine 
funktionentheorie und elliptische funktionen}, Springer-Verlag, 1964

\bibitem{H1} Y.~Haraoka, {\it Canonical forms of differential equations free from
accessory parameters}, SIAM J. of Mathematical Analysis, vol.25, n.4 (1994), 1203-1226

\bibitem{H2} Y.~Haraoka, {\it Monodromy representations of systems of differential 
equations free from accessory parameters}, SIAM J. Math. Anal. 25 (6) (1994) 
1595-1621

\bibitem{HY1} Y.~Haraoka and T.~Yokoyama, {\it Construction of rigid local systems
and integral representation of their sections}, to appear in Math.$\!$ Nachr.

\bibitem{Hit1} N.~Hitchin, {\it Hyperk\"ahler manifolds}, s\'eminaire Bourbaki,
44\'eme ann\'ee, 1991-92, n.~748

\bibitem{NKatz} N.~Katz, {\it Rigid local systems}, Annals of Mathematics Studies, n.139, 
Princeton University Press (1996).

\bibitem{Ko1} V.P.~Kostov, {\it Monodromy groups of regular systems on Riemann's
sphere}, prepublication N 401 de l'Universite de Nice, 1994; to appear in
Encyclopedia of Mathematical Sciences, Springer

\bibitem{Ko2} V.P.~Kostov, {\it The Deligne-Simpson problem -- a survey},
J.$\!$ of Algebra, in press; available online at {\it www.sciencedirect.com}

\bibitem{Mag1} P.~Magyar, {\it Bruhat order for two flags and a line}, Journal of
Algebraic Combinatorics 21 (2005)

\bibitem{MWZ} P.~Magyar, J.~Weyman, A.~Zelevinsky, {\it Multiple flag varieties of
finite type}, Advances in Mathematics 141, 97-118 (1999)

% \bibitem{N1} H.~Nakajima, {\it Instantons on ALE spaces, quiver varieties, and 
% Kac-Moody algebras}, Duke math. journal, vol. 76, n. 2, 1994 

% \bibitem{N2} H.~Nakajima, {\it Quiver varieties, and Kac-Moody algebras}, Duke math. journal, 
% vol. 91, n. 3, 1998

\bibitem{O} K.~Okubo, {\it On the group of Fuchsian equations}, Seminar Reports
of Tokyo Metropolitan University, 1987

\bibitem{OrSch1} A.~Orlov, D.~Scherbin, {\it Multivariate hypergeometric functions
as tau functions of Toda lattice and Kadomtsev-Petviashvili equation}, preprint
arXiv:math-ph/0003011 v2, Oct. 19, 2000

\bibitem{OrSch2} A.~Orlov, D.~Scherbin, {\it Fermionic representations for basic
hypergeometric functions related to Schur polynomials}, preprint
arXiv:nlin.SI/0001011 v4, Dec. 1, 2000

\bibitem{SJM} M.~Sato, M.~Jimbo, T.~Miwa, {\it Holonomic quantum fields I, III, IV, V}
(in Russian), Mir publishing house, Moscow, Russia, 1983 

\bibitem{S} C.~Simpson, {\it Products of matrices}, AMS Proceedings 1 (1992),
157-185

% \bibitem{V} V.S.~Varadarajan, {\it Meromorphic differential equations}, 
% Expositiones Mathematicae, International Journal for Pure and Applied
% Mathematics, Vol. 9, Number 2, 1991

\bibitem{WW} E.~Whittaker, G.~Watson, {\it A Course of modern analysis},
Cambridge University Press, 1927

\bibitem{Wil} G.~Wilson, {\it Collision of Calogero-Moser particles and an adelic 
Grassmanian}, Invent. math. 133, 1-41 (1998)

\end{thebibliography}
\end{document}